\documentclass[11pt]{amsart}
\usepackage{hyphenat}
\usepackage[kerning, tracking, spacing]{microtype}
\usepackage[margin=1in]{geometry}
\usepackage{appendix}
\usepackage{amsmath, amsthm}
\usepackage{amssymb}
\usepackage{mathrsfs}
\usepackage{enumitem}
\usepackage{graphicx}
\usepackage{mathdots}
\usepackage{color}
\usepackage{hyperref,url}
\usepackage{tikz, tikz-cd}

\usetikzlibrary{matrix}
\usetikzlibrary{shapes}
\usetikzlibrary{arrows,decorations.markings, tikzmark}
\usepackage{mathtools}
\usepackage{ stmaryrd }
\usepackage[normalem]{ulem}
\usepackage[utf8]{inputenc}

\usepackage[notcite, notref]{showkeys}

\usepackage{xcolor}
\hypersetup{
    colorlinks,
    linkcolor={red!50!black},
    citecolor={blue!50!black},
    urlcolor={blue!80!black}
}

\def\Mc{\mathcal{M}}
\def\CH{\mathsf{CH}}

\newcommand{\Sym}{\operatorname{Sym}}

\def\Ac{\mathcal{A}}

\newcommand{\taut}{\operatorname{taut}}
\newcommand{\Tor}{\operatorname{Tor}}

\def\CC{\mathbb{C}}

\def\EE{\mathbb{E}}

\def\RR{\mathbb{R}}
\def\QQ{\mathbb{Q}}

\def\ZZ{\mathbb{Z}}

\def\SL{\operatorname{SL}}

\def\Hom{\operatorname{Hom}}

\def\Abar{\overline{\mathcal{A}}}

\newcommand{\bbC}{\mathbb C}
\newcommand{\bbZ}{\mathbb Z}

\newcommand{\NL}{\mathsf{NL}}
\newcommand{\KM}{\mathrm{KM}}

\newcommand{\Sp}{\operatorname{Sp}}

\renewcommand{\d}{\mathrm{d}}

\newcommand{\ddc}{\mathrm{d}\,\mathrm{d}^c\, }

\newcommand{\str}{\operatorname{tr}_s}

\newcommand{\ii}{\iota}
\newcommand{\eps}{\varepsilon}

\usepackage{amsthm}
\theoremstyle{definition}
\newtheorem{definition}{Definition}[section]

\newtheorem{theorem}[definition]{Theorem}

\newtheorem{proposition}[definition]{Proposition}
\newtheorem{corollary}[definition]{Corollary}
\newtheorem{lemma}[definition]{Lemma}
\newtheorem{conjecture}{Conjecture}

\newtheorem{remark}[definition]{Remark}
\newtheorem{notation}[definition]{Notation}

\newtheorem*{conjecture*}{Conjecture}
\newtheorem*{question*}{Question}
\newtheorem*{theorem*}{Theorem}

\newtheorem{introthm}{Theorem}

\newtheorem{introcor}{Corollary}[introthm]

\title{Mathai--Quillen on locally symmetric spaces}
\author{}
\date{}

\microtypecontext{spacing=nonfrench}

\title{Symplectic Noether--Lefschetz modularity}
\author{Luis Garc\'{i}a, Fran\c{c}ois Greer, Aitor Iribar L\'{o}pez}
\date{}

\begin{document}

\begin{abstract}
    We study the locus $\NL_{g,d}$ of principally polarized abelian varieties of dimension $g$ admitting a degree $d$ homomorphism from an elliptic curve. We show that the associated cycle class $[\NL_{g,d}]$ is the $d$-th Fourier coefficient of a modular form of weight $2g$ with values in the cohomology of $\mathcal A_g$. If $g =2,3,4,5,6,7,12$ or any even number $16 \leq g \leq 6000$, then modularity is the only source of relations between the Noether--Lefschetz cycles.
\end{abstract}

\maketitle

\setcounter{tocdepth}{1}
\tableofcontents

\section{Introduction}

\subsection{Abelian varieties and modular forms}
Let $\Ac_g$ be the moduli stack of complex principally polarized abelian varieties of dimension $g \geq 1$. Any line bundle on an abelian variety $A$ induces a map from $A$ to its dual $A^\vee$. A principal polarization on $A$ is the Chern class $\theta_A = c_1(\mathcal L) \in H^2(A, \ZZ)$ of an ample line bundle $\mathcal L$ that induces an isomorphism from $A$ to its dual. The data of $\theta_A$ is equivalent to the data of the map $A \to A^\vee$, and we will denote both of them by $\theta_A$. Complex points of $\Ac_g$ are therefore in bijection with isomorphism classes of such pairs $(A, \theta_A)$.

Over the complex numbers, $\Ac_g$ can be realized as the orbifold quotient
$$
[\mathfrak H_g / \Sp_{2g}(\ZZ)]\,.
$$
Here, $\mathfrak H_g$ is the Siegel space of complex, symmetric $g \times g$ matrices $Z = X+iY$, where $Y$ is positive definite, $\Sp_{2g}(\ZZ)$ is the group of $2g\times 2g$ symplectic matrices with integer coefficients, and the action is given by
$$
\begin{pmatrix}
    A&B \\C&D
\end{pmatrix}\cdot Z = (AZ+B)(CZ+D)^{-1}\,, \quad \begin{pmatrix}
    A&B \\C&D
\end{pmatrix} \in \Sp_{2g}(\ZZ)\,.
$$

The moduli stacks $\Ac_g$ are not proper, but they admit several compactifications (\cite{bb}, \cite{AMRT}, \cite{FaltingsChai}). In this article, we will only encounter compactifications of $\Ac_1$. In this case, all of the above compactifications agree, and the resulting smooth, proper stack is denoted by $\Abar_1$. Its boundary $\Abar_1 \smallsetminus \Ac_1$ is a single point, which we will denote by $i\infty$.

The Hodge bundle $\mathbb E_g \to \Ac_g$ is a vector bundle of rank $g$ with the property
$$
\mathbb E_g|_{(A,\theta_A)} = H^0(A, \Omega_A).
$$
It plays a major role in the study of the $\Ac_g$. It extends naturally to any toroidal compactification; in particular, it extends to $\Abar_1$. Sections of $\EE_1^{\otimes k}$ on $\Abar_1$ correspond to modular forms of weight $k$ for the group $\Sp_{2}(\ZZ) = \SL_2(\ZZ)$; that is, holomorphic functions $f : \mathfrak H_1 \to \CC$ such that
$$
f \left( \frac{a\tau + b}{c\tau +d}\right) = (c\tau + d)^{k}f(\tau)
$$
for all $a,b,c,d \in \ZZ$ with $ad-bc = 1$, and such that $f$ is bounded as $y = \operatorname{Im}(\tau) \to \infty$. 

A modular form is periodic, so it has a Fourier expansion, which can be viewed as a Taylor expansion at the cusp $i\infty$:
$$
f (\tau) = \sum_{d =0}^\infty a_d q^d\,,\quad q = e^{2\pi i \tau}
$$
and we say that $f$ is rational if all of its Fourier coefficients are rational. Let $\operatorname{Mod}_{k}$ be the space of rational modular forms of weight $k$. Since $\mathcal A_1$ and $\EE_1$ are defined over $\QQ$, every modular form is a complex linear combination of rational modular forms. A cusp form is a modular form whose constant Fourier coefficient $a_0$ vanishes. The vector space of cusp forms of weight $k$ is denoted by $\operatorname{SMod}_k$.

More generally, sections of the Hodge bundle $\EE_g \to \Ac_g$ or other vector bundles obtained from it (e.g. the determinant) are called \emph{Siegel modular forms} of genus $g$.

\subsection{The cohomology of $\Ac_g$}

The Chern classes $\lambda_i = c_i(\mathbb E)$, $i =1,\ldots , g$, generate subrings of the Chow and cohomology rings\footnote{In this paper, Chow and cohomology rings will be taken with rational or even complex coefficients. In the latter case, we will identify them with the de Rham cohomology without explicit mention.} of $\Ac_g$, which are usually referred to as the \emph{tautological rings}:
$$
R^*(\Ac_g) \subset \CH^*(\Ac_g)\,, \quad RH^*(\Ac_g) \subset H^*(\Ac_g)\,.
$$
These tautological classes form the simplest collection of algebraic cycles on $\Ac_g$.

By Borel's stability theorem \cite{Borel}, every cohomology class of small dimension is tautological:
$$
H^k(\Ac_g) = RH^k(\Ac_g)\,, \quad k < g.
$$
An improvement by Ta\"ibi \cite[Theorem 33]{HT18} shows that this bound can be substantially improved if we restrict to the \emph{pure} cohomology of $\Ac_g$: $H^k_{\mathsf{pure}}(\Ac_g)$ is the subspace of $H^k(\Ac_g)$ generated by restrictions of cohomology classes from any toroidal compactification $\Abar_g^{\Sigma}$. It is a sub-Hodge structure of pure weight $k$, and it contains the cycle associated to any algebraic subvariety of $\Ac_g$. Ta\"ibi proves that
\begin{equation}\label{eq:taibi-bound}
H^k_{\mathsf{pure}}(\Ac_g) = RH^k(\Ac_g), \quad  k <2g-2.
\end{equation}
In this article, we study another collection of algebraic cycles on $\Ac_g$, coming from abelian varieties containing an elliptic curve. The simplest of these cycles is given by the \emph{product locus}
$$
i: \Ac_1 \times \Ac_{g-1} \longrightarrow \Ac_g\,, \quad [\Ac_1 \times \Ac_{g-1}] = i_*(1)\in \CH^{g-1}(\Ac_g) \to H^{2g-2}(\Ac_{g})
$$
By \cite{COP24} and \cite{IribarLopez}, the cycle $[\Ac_1 \times \Ac_{g-1}]$ is \emph{not} tautological if $g=6, 12$ or $g \geq 16$ even.

More generally, we will consider the cycle $\NL_{g,d}$, parametrizing quadruples $(A, \theta_A, E, f)$, where $(A, \theta_A)$ is a principally polarized abelian variety of dimension $g$, $E$ is an elliptic curve, $f :E \to A$ is a homomorphism, and the composition $f^\vee \circ \theta_A \circ f$ is multiplication by $d$ on $E$, under the canonical identification $E \cong E^\vee$.

The map $i : \NL_{g,d} \to \mathcal A_g$ that remembers $(A, \theta_A)$ is proper, of codimension $g-1$. Its image is part of the \emph{Noether--Lefschetz locus} on $\Ac_g$, which parametrizes abelian varieties whose Picard rank is at least $2$. Let
$$
[\NL_{g,d}] = i_*(1) \in \CH^{g-1}(\Ac_g) \to H^{2g-2}(\Ac_g)
$$
be the associated \emph{Noether--Lefschetz} cycle class. We can arrange these classes into the \emph{Noether--Lefschetz} power series $\Theta_{\NL}$:
$$
\Theta_{\NL} (q) = \frac{(-1)^g}{24}\lambda_{g-1} + \sum_{d \geq 1} [\NL_{g,d}]q^d \in \CH^{g-1}(\Ac_g) \llbracket q \rrbracket\,.
$$
It will be convenient to define $[\NL_{g,0}]=(-1)^g \lambda_{g-1}/24$. The following conjecture appears in \cite{GreerLian, IribarLopez}, and was also suggested by A. Pixton.

\begin{conjecture}\label{conjecture:modularity}
    The Noether--Lefschetz series $\Theta_{\NL}(q)$ is the Fourier expansion of a modular form of weight $2g$ for $\SL_2(\ZZ)$ with values\footnote{If $V$ is any vector space, an element $f \in V\llbracket q \rrbracket$ is a modular form of weight $k$ with values in $V$ if for every linear functional $\mu \in V^\vee$, $\mu(f)$ is the Fourier expansion of a modular form of weight $k$.} in $\CH^{g-1}(\Ac_g)$.
\end{conjecture}

Our main result in this article is the proof of Conjecture \ref{conjecture:modularity} in cohomology:
\begin{introthm}\label{thm:main}
    The image of $\Theta_{\NL}$ under the cycle class map
    $$
    \CH^{g-1}(\Ac_g) \to H^{2g-2}(\Ac_{g})
    $$
    is the Fourier expansion of a modular form of weight $2g$ for $\SL_2(\ZZ)$, for $g \geq 2$.
\end{introthm}

During the preparation of this article, T. Cheek, R. K. Tamura, and J. Xia informed us that they found a proof of Theorem \ref{thm:main} independently, using similar ideas from Kudla--Millson theory and its reformulation in terms of Mathai--Quillen forms from \cite{Garciasuper,  Branchereau}; see \cite{TimRyan}. Subsequently, we found out that H. Yao gave another proof of Theorem \ref{thm:main} using the Borel--Serre compactification of $\mathcal A_g$; see \cite{Hanyu26}.

Let $\mathbf{NL}_g \subset H^{2g-2}(\mathcal A_g)$ be the $\QQ$-vector space generated by the cycles $[\NL_{g,d}]$, $d \geq 0$. Theorem \ref{thm:main} can be reformulated as saying that there is a well-defined linear map
\begin{equation}\label{eq:Psi-map}
     \psi: (\operatorname{Mod}_{2g})^\vee \longrightarrow \mathbf{NL}_g\,, \quad a_d(-) \mapsto [\NL_{g,d}]\,,
\end{equation}
where $a_d(-) : \operatorname{Mod}_{2g} \to \QQ$ is the linear map that extracts the $d$-th Fourier coefficient. Clearly, $\psi$ is surjective by definition, so we have
\begin{equation}\label{eq:dim-bound}
\dim(\mathbf{NL}_g) \leq \dim(\operatorname{Mod}_{2g})\,.
\end{equation}
Ta\"ibi's calculation in \cite[Proof of Theorem 33]{HT18} and the
calculation for parameters of type~C in
\cite[proof of Theorem~2.4]{CanningPetersenTaibi2026} show that
$$
\dim\left(\frac{H_{\mathsf{pure}}^{2g-2}(\Ac_g)}{RH^{2g-2}(\Ac_g)}\right) \leq \dim(\operatorname{SMod}_{2g})\,.
$$
Therefore, if \eqref{eq:dim-bound} is an equality, the pure cohomology group $H^{2g-2}_{\mathsf{pure}}(\Ac_g)$ is generated by tautological classes and Noether--Lefschetz cycles $[\NL_{g,d}]$; see Remark \ref{rem:NL-taut}. This would explain the pure cohomology of $\Ac_g$ in the first case beyond \eqref{eq:taibi-bound}. We conjecture that this is the case:
\begin{conjecture}\label{conjecture-NL}
    For all $g \geq 2$, we have
    $$
    \dim(\mathbf{NL}_g) = \dim(\operatorname{Mod}_{2g}) = \begin{cases}
    \lfloor g/6 \rfloor & \text{if } g \equiv 1 \pmod 6\,,\\
    \lfloor g/6 \rfloor + 1 & \text{otherwise}\,.
\end{cases}
    $$
    In particular, the pure cohomology group $H^{2g-2}_{\mathsf{pure}}(\Ac_g)$ is generated by tautological classes and the cycles $[\NL_{g,d}]$.
\end{conjecture}
\noindent Note that a similar conjecture that explains the low degree cohomology of orthogonal Shimura varieties in terms of special cycles has been proven in \cite{BergeronLiMM}.

\subsection{Hecke operators}

Let $\mathcal H (\SL_2(\ZZ))$ be the classical algebra of Hecke operators, spanned by operators $T_n$, $n \in \ZZ_{>0}$. It acts on the space of modular forms $\operatorname{Mod}_{2g}$, and it preserves the space of cusp forms $\operatorname{SMod}_{2g}$.

By the multiplicity one theorem for modular forms \cite{diamondshurman}, the action of $\mathcal H (\SL_2(\ZZ))$ is semisimple, so $\mathrm{SMod}_{2g}\otimes \CC$ decomposes into $1$-dimensional eigenspaces for the action of $\mathcal H(\SL_{2}(\ZZ))$. Over the rationals, however, this is not true. In fact, the following has been conjectured by Maeda \cite{HM97}:
\begin{conjecture}[Maeda]\label{conj:modular-forms}
    As an algebra over $\mathcal H (\SL_2(\ZZ))$, the space of cusp forms $\operatorname{SMod}_{2g}$ is simple.
\end{conjecture}
Conjecture \ref{conj:modular-forms} has been verified in \cite{GM12} for $2g \leq 12000$ by checking that the characteristic polynomial of $T_2$ is irreducible over $\QQ$. In Section \ref{sec:hecke} we find a compatible geometric Hecke algebra action on $\mathbf{NL}_g$ that makes $\psi$ equivariant. As a consequence, we obtain:
\begin{introthm}\label{thm:non-vanish}
    Let $g \geq 2$. Suppose that Conjecture \ref{conj:modular-forms} holds for $\operatorname{SMod}_{2g}$. Then one of the following conditions must hold:
    \begin{itemize}
        \item[(a)] Conjecture \ref{conjecture-NL} holds on $\Ac_g$, or
        \item[(b)]  $\mathbf{NL}_g = \QQ\lambda_{g-1}$ is $1$-dimensional.
    \end{itemize}
\end{introthm}
As evidence for Conjecture \ref{conjecture-NL}, we use Theorem \ref{thm:non-vanish} to prove it in a number of cases:
\begin{introcor}\label{cor:psi-is-iso}
    Conjecture \ref{conjecture-NL} is true if $g = 2,3,4,5,6,7,12$ or $16 \leq g \leq 6000$ even.
\end{introcor}

\subsection{Effective representatives of $\lambda_{g-1}$}

If $g =2,3$, it has been shown that $\lambda_{g-1}$ is a positive multiple of the class $[\mathcal A_1 \times \mathcal A_{g-1}] = [\NL_{g,1}]$:
$$
\lambda_1 = \frac{1}{10}[\mathcal A_1 \times \Ac_1]\quad\text{on } \Ac_2\,, \quad \lambda_2 = \frac{1}{21}[\Ac_1 \times \Ac_2]\quad\text{on } \Ac_3\,.
$$
By Theorem \ref{thm:main}, we also have
\begin{equation}\label{eq:low-g}
\lambda_{g-1} = \frac{6|B_{2g}|}{g\sigma_{2g-1}(d)}[\NL_{g,d}]
\end{equation}
for all $d \geq 1$ whenever the space of modular forms of weight $2g$ is $1$-dimensional, i.e., for $g = 2,3,4,5,7$. This also follows from Ta\"ibi's work and the tautological projection calculations of \cite{IribarLopez}.

Note that the cycles $\NL_{g,d}$ are reducible, as noted in \cite{GreerLian, IribarLopez}. If we let $\NL_{g,d}^{\mathsf{prim}}$ parametrize principally polarized abelian varieties $(A, \theta_A)$ containing an elliptic subgroup $E \subset A$ such that $\int_A[E]\cdot \theta_A=d$, then $\NL_{g,d}^{\mathsf{prim}}$ is an irreducible subvariety of $\Ac_g$, of codimension $g-1$. As cycles, we have
\begin{equation}\label{eq:NL-from-prim}
    [\NL_{g,d}] = \sum_{e \mid d} \sigma_1(d/e) [\NL_{g,e}^{\mathsf{prim}}]\,.
\end{equation}
Following \cite{BM19}, we study the cone spanned by the primitive $\NL$ cycles:
\begin{introthm}\label{thm:lg-1}
    Let $C, C^{\mathsf{prim}} \subset H^{2g-2}(\mathcal A_g, \mathbb R
    )$ be the cones spanned by the cycles $[\NL_{g,d}]$ and $[\NL_{g,d}^{\mathsf{prim}}]$ respectively, for $d\geq 1$. Then $C$ and $C^{\mathsf{prim}}$ are rational polyhedral cones, and $\lambda_{g-1}$ lies in their interior.
\end{introthm}

\subsection{$d$-elliptic curves}

Let $\Mc_g$ denote the moduli space of genus $g$ curves. The analogue of the space $\NL_{g,d}$ for curves is the \emph{Hurwitz space}
$$
\mathcal H_{g\to 1, d}\,.
$$
It parametrizes degree $d$ maps $C \to E$ from a genus $g$ curve to an elliptic curve $E$. It is a smooth, equidimensional orbifold. The forgetful map
$$
\mathsf{ft}: \mathcal H_{g \to 1, d} \to \mathcal M_g
$$
is proper of relative codimension $g-1$. Its image is the locus of of $d$-elliptic curves. Let
$$
[\mathcal H_{g \to 1, d}] = \mathsf{ft}_*(1) \in H^{2g-2}(\mathcal M_g)
$$
be the associated cycle. As a consequence of Theorem \ref{thm:main}, we obtain:
\begin{introthm}\label{thm:curves}
    For $g \geq 2$, the generating series
    $$
    \sum_{d \geq 2} [\mathcal H_{g \to 1, d}] q^d \in H^{2g-2}(\mathcal M_g)\llbracket q \rrbracket
    $$
    is the Fourier expansion of a cusp modular form of weight $2g$ for $\SL_2(\ZZ)$.
\end{introthm}

The moduli space of curves also has a natural \emph{tautological ring} $R^*(\Mc_g)$, and the locus of bielliptic curves (i.e. $2$-elliptic curves) has been a source of \emph{non-tautological} classes. We prove:
\begin{introcor}\label{cor:curve-vanish}
    If $g \leq 11$, or $g =13$, then the fundamental class $[\mathcal H_{g\to 1, d}] \in H^{2g-2}(\mathcal M_g)$ vanishes. 
\end{introcor}

In genus $12$, consider the square of the discriminant form
$$
\Delta(q)^2 =q^2\prod_{n \geq 1}(1-q^n)^{48} =  \sum_{d \geq 2} a_d(\Delta^2) q^d
$$
It is conjectured that $a_d(\Delta^2)\neq 0$. In fact, $a_d\neq 0$ if and only if the characteristic polynomial of $T_d$ acting on $\operatorname{SMod}_{24}$ is irreducible over $\QQ$.
The square of the discriminant had appeared in work of C. Lian \cite{Lian22}, where he reduced the failure of the closure $\overline{\mathcal H}_{g\to 1, d}\subset \overline{\mathcal M}_g$ of our Hurwitz cycle to be tautological in the Deligne--Mumford compactification to the nonvanishing of $a_d(\Delta^2)$, for $g \geq 12$. We obtain the following related result:
\begin{introcor}\label{cor:Lian-smooth}
    For any $d \geq 2$, the class $[\mathcal H_{12\to 1,d}] \in H^{22}(\mathcal M_{12})$ is nonzero if and only if $a_d(\Delta^2) \neq 0$.
\end{introcor}

\subsection{Organization of the paper} The proof of the main result uses the Kudla--Millson form on $\mathcal A_1 \times \mathcal A_g$, as suggested in \cite{GreerLian}. We show that its salient properties extend across the cusp of $\mathcal A_1$. It is convenient for the proof of this extension to formulate the Kudla--Millson form in the Mathai--Quillen formalism from \cite{Garcia}. In \S \ref{sec:notation} we introduce the Mathai--Quillen forms that we will use in the proof of Theorem \ref{thm:main}, and write them down explicitly on the universal cover of $\mathcal A_1 \times \mathcal A_g$. In \S \ref{sec:proof}, we perform cusp estimates to prove Theorem \ref{thm:main}. In \S \ref{sec:hecke}, we introduce basic Hecke operators and prove Theorem \ref{thm:non-vanish} and Corollary \ref{cor:psi-is-iso}. Finally, \S \ref{sec:geometry} is dedicated to the proof of  Theorems \ref{thm:lg-1} and \ref{thm:curves}.

\subsection{Acknowledgements}
We thank the organizers of the conference \emph{Moduli Spaces and Modular Forms} (2024) on Schiermonnikoog island, where conversations about using the methods of \cite{Garcia} to prove modularity on $\Ac_g$ started. We also thank I. Barros, S. Canning, J. Feusi, C. Lian, S. Tayou, and G. van der Geer for enlightening conversations about the cohomology of $\Ac_g$ and other Shimura varieties.

A.I.L. would like to thank R. Pandharipande for suggesting to work on $\NL$ cycles on $\Ac_g$ as part of his doctoral studies at ETH Z\"{u}rich. A.I.L. was supported by SNF Mobility, grant number P500--2$\_$239077.

We are grateful to T. Cheek, R. K. Tamura and J. Xia for conversations around the topic of this work. We thank them, as well as H. Yao, for coordinating the release of the preprints.

\subsection{AI Disclosure} 
The approaches to all the theorems in the paper come from the authors, and we have written all of the proofs as we understand them. The evaluation in Proposition \ref{p:W2exact0} was first found using A.I.

\section{Notation and first remarks}\label{sec:notation}

To avoid confusion, we will use $\mathfrak H_1$ to denote the Siegel upper half-space that gives rise to $\Ac_1, \Abar_1$, and $\mathscr{H}$ for the space for the modular parameter $\tau = u+iv$.

\subsection{Mathai--Quillen forms}

Assume we have:
\begin{itemize}
    \item A rank $r$ holomorphic vector bundle $\mathbb F$ on a complex manifold $M$. 
    \item A Hermitian metric $h$ on $\mathbb F$. The dual bundle $\mathbb F^\vee$ inherits the dual metric: if (choosing a trivialization of $\mathbb F$) the metric on $\mathbb F$ is given by $\|s\|^2 = s^* hs$ with $h$ a positive-definite Hermitian $r \times r$ matrix, then the dual metric is given by $\|\alpha\|^2=\alpha^* h^{-1}\alpha$. The associated Chern--Weil connection is 
    $$
    \nabla_{\mathbb F^\vee} = \nabla^{1,0}+\nabla^{0,1} : \mathbb F^\vee \otimes A^*(M) \longrightarrow \mathbb F^\vee\otimes A^{*+1}(M)
    $$
    with curvature $\nabla^2_{\mathbb F^\vee}$. We will omit the subscript indices when they are clear from the context.
    \item A holomorphic section $s :   \mathcal O_M \to \mathbb F$, with adjoint $s^{\dagger} :  \mathbb F \to \mathcal O$. Let $\|s\|^2 = s^\dagger s$ be the square norm with respect to $h$.
\end{itemize}
From this setup, Mathai and Quillen construct in \cite{Quillen1985, MathaiQuillen1986} a superconnection on the superalgebra $\wedge^*\mathbb F^\vee \otimes A^*(M)$ by the formula
$$
\nabla_{s} = \nabla_{\wedge^*\mathbb F^\vee} + \sqrt{-2\pi}(\iota_{s} + s^\dagger \wedge -)\,.
$$
See also \cite{Garciasuper}. The modified Euler class of $\mathbb F$ is defined to be
\begin{equation}\label{eq:define-mod-euler}
\mathsf{e} (\mathbb F, \nabla, s) = \mathsf{e}(\mathbb F, s) =  \left(\frac{-1}{2\pi i}\right)^{r}\left[\str \left(e^{\nabla^2_s}\right)\right]_{2r} \in A^{r,r}(M)\, .
\end{equation}
\begin{lemma}\label{lem:ex-euler}
    Let $B_s$ be the $r \times r$ matrix of $(1,1)$ forms given by
    $$
    B_s = \nabla^{1,0}_{\mathbb F}s\wedge \nabla_{\mathbb F^\vee}^{0,1}s^{\dagger} \in A^{1,1}(\mathbb F \otimes \mathbb F^\vee) = A^{1,1}(\operatorname{End}(\mathbb F))\, ,
    $$
    and let $R$ be the $r \times r$ matrix of the curvature $\nabla_{\mathbb F}^2\in A^{1,1}(\operatorname{End}(\mathbb F))$.
    Then,
    $$
    \mathsf{e}(\mathbb F, s) = \left(\frac{-1}{2\pi i}\right)^re^{-2\pi\|s\|^2}\sum_{k =0}^r \frac{(2\pi)^k}{k!}[t^k]\det(R + tB_s)\, ,
    $$
    where $[t^k]$ denotes extraction of the coefficient of $t^k$.
\end{lemma}
\begin{proof}
    In a local frame given by $e_1, \ldots , e_r$, the superconnection is given as
    \begin{equation}\label{eq:local-superconnection}
    \nabla^2_{s} = \sum R_{ij} \eps_j \ii_i + \sqrt{-2\pi}\left(\sum a_i\ii_i + b_i \eps_i\right) - 2\pi\|s\|^2 \operatorname{id}\,,
    \end{equation}
    where $\ii_i = \iota_{e_i^\vee}$, $\eps_i = e_i \wedge -$, $a = \nabla^{1,0}s$ and $b = \nabla^{0,1}s^\dagger$. They satisfy
    \begin{equation}\label{eq:Clifford}
        \eps_i \eps_j + \eps_j \eps_i =\ii_i\ii_j+\ii_j\ii_i=0\, \quad \ii_i\eps_j+\eps_j\ii_i = \delta_{ij}\,.
    \end{equation}
    Note that all the terms in \eqref{eq:local-superconnection} have the same $A^*$-degree as the degree of the monomial in the $\iota_i$ and $\eps_i$. Since the Euler class \eqref{eq:define-mod-euler} extracts the $A^*$-degree $2r$, it follows that, to calculate the supertrace of the expansion of  $\exp(\nabla_s^2)$, we can restrict to terms that have words in the $\eps_i$ and $\iota_j$ of length $\leq 2r$. Using the relations \eqref{eq:Clifford}, any term can be uniquely written in terms of the ordered monomials
    $$
    M_{I,J} = \eps_I \ii_J\, , \quad \eps_I=\eps_{i_1}\cdots \eps_{i_k}\,, \qquad\ii_I=\ii_{i_k}\cdots \ii_{i_1}\,.
    $$
    Since the supertrace is $A^*$-linear, it is enough to compute it on these endomorphisms. An explicit computation shows that
    $$
    \str(M_{I,J})=\begin{cases}(-1)^r, & I=J=\{1,\dots,r\},\\
    0, & \text{otherwise}\,.\end{cases}
    $$
    Therefore, since the term $\delta_{ij}$ in \eqref{eq:Clifford} decreases the degree in $\eps_i$ and $\iota_j$, it follows that, for the purposes of the calculation of $\mathsf{e}(\mathbb F,s)$, we can assume that all the $\eps_i$ and $\iota_j$ supercommute with each other:
    $$
    \iota_i \eps_j + \eps_j \iota_i = 0\,.
    $$
    Thus, the desired supertrace term is
    $$
    \mathsf{e}^{\circ}(\mathbb F, s) = e^{-2\pi\|s\|^2}\str(M_{[r], [r]})\sum_{\ell=0}^{r}[\eps_{[r]} \ii_{[r]}]\frac{\left(\sum_{i,j}R_{ij}\eps_i \ii_j\right)^{r-\ell}}{(r-\ell)!}\frac{(\sqrt{-2\pi}\sum_{i}a_i \ii_i)^\ell}{\ell!}\frac{(\sqrt{-2\pi}\sum_{i}b_i \eps_i)^\ell}{\ell!}\,.
    $$
    We have
    $$
    \Big(\sum_i a_i \ii_i\big)\Big(\sum_i b_j \eps_j\big) = \sum_{i,j}a_i \wedge b_j\eps_i \ii_j\, ,
    $$
    and so
    $$
    \mathsf{e}^{\circ}(\mathbb F, s) = e^{-2\pi \|s\|^2}(-1)^r\sum_{\ell=0}^{r}\frac{(-2\pi)^\ell}{\ell!}[\eps_{[r]} \ii_{[r]}]\frac{\left(\sum_{i,j}R_{ij}\eps_i \ii_j\right)^{r-\ell}}{(r-\ell)!}\frac{(\sum_{i,j}a_i\wedge b_j \eps_i\ii_j)^\ell}{\ell!}\,,
    $$
    which is the determinant formula from the statement, since
    $$
    \frac{1}{r!}\left(\sum_{i,j}m_{ij}\eps_i \ii_j\right)^r = \det(M)\eps_{[r]}\ii_{[r]}\,.
    $$
\end{proof}

\begin{definition}
    Let $\lambda$ be a formal parameter. The $\lambda$-deformation of the modified Euler class is
    $$
    \mathsf{e}_\lambda(\mathbb F, \nabla, s) = \left(\frac{-1}{2\pi i }\right)^re^{-2\pi\|s\|^2}\sum_{\ell =0}^r \frac{(2\pi)^\ell}{\ell!}[t^\ell]\det(R -2\pi i \lambda \operatorname{Id} + tB_s ) \in A^{*,*}(M)[\lambda]\,,
    $$
    where $R$ and $B_s$ are as in Lemma \ref{lem:ex-euler}. It is a homogeneous $(r,r)$-form, if we think of $\lambda$ as having bidegree $(1,1)$.
\end{definition}

In the Mathai--Quillen formalism, the double transgression for the modified Euler class is usually given by the form
$$
\left(\frac{-1}{2\pi i}\right)^{r-1}\left[\str\left(\mathcal Ne^{\nabla^2_s}\right)\right]_{2r-2}\in A^{r-1,r-1}(M)\,,
$$
where $\mathcal N$ is the operator acting by multiplication by $-k$ on $\wedge^k \mathbb F^\vee$. We will write this form differently:
\begin{lemma}\label{lem:ex-green}
    We have
    $$
    \left(\frac{-1}{2\pi i}\right)^{r-1}\left[\str\left(\mathcal Ne^{\nabla^2_s}\right)\right]_{2r-2} = \operatorname{Coeff}(\lambda^1; \mathsf{e} _\lambda(\mathbb F,\nabla,s))
    $$
\end{lemma}
\begin{proof}
    The proof is the same as in Lemma \ref{lem:ex-euler}.
\end{proof}

\begin{definition}\label{def:Green}
    The Green function is defined as
$$
\mathfrak g(\mathbb F, s) := [\lambda^1]\int_{1}^{\infty}\mathsf{e}_{\lambda}(\mathbb F, t^{1/2}s)\frac{dt}{t}\, .
$$
\end{definition}

Since $\det(\lambda \mathrm{Id} +  N) = \det(M) + \lambda \mathrm{tr}(\mathrm{adj} (M)) + \ldots$, we have
\begin{equation}\label{eq:green-explicit}
\mathfrak g(\mathbb F, s) =\left(\frac{-1}{2\pi i}\right)^{r-1}\sum_{k=0}^{r-1} \frac{2\pi}{k!}[t^k]\operatorname{tr}(\operatorname{adj}(R + tB_s))\int_{1}^{\infty}e^{-2\pi \sigma\|s\|^2}(2\pi \sigma)^{k-1}\, d\sigma\,,
\end{equation}
which converges whenever $s\neq 0$.

\begin{lemma}\label{lem:green-is-current}
    If $s$ is a regular section\footnote{Here, regular section should mean that the zero section and $s$ are transverse at the level of tangent spaces.}, then $\mathfrak g(\mathbb F, s)$ is a current on $M$.
\end{lemma}
\begin{proof}
Since $s$ is a regular section, by the implicit function theorem the local model is given by
$$
C_K\int_1^{\infty}e^{-t\|z_1\|^2}t^{\ell-1}dt\,,\quad  0 \leq\ell \leq r-1\, , z=(z_1, z_2) \in \mathbb C^r \times K \subset \mathbb C^r\times \mathbb C^{\dim M - r}\, , K\text{ compact}
$$
and we have, by the change $u = 2\pi t\|z_1\|^2$
\begin{align*}\label{eq:inc-Gamma}
2\pi \int_1^{\infty}e^{-2\pi t\|z_1\|^2}(2\pi t)^{\ell-1}dt =\,  \|z_1\|^{-2\ell}\int_{2\pi \|z_1\|^2}^{\infty}e^{-u}u^{\ell-1}du 
\longrightarrow \begin{cases}
    \|z_1\|^{-2\ell}(\ell-1)!&\text{ if }\ell >0\\
    \log \|z_1\|^{-2}+O(1)&\text{ if }\ell=0
\end{cases}\,.
\end{align*}
Both expressions are locally integrable on $\mathbb C^r$, as long as $\ell <r$.
\end{proof}
\begin{lemma}\label{l:g-curr-dZ}
    Let $\mathrm{d}^c= \frac{1}{4\pi i}(\partial -\overline{\partial})$. If $s$ is regular around the zero section then, as currents,
    $$
    \ddc \mathfrak g(\mathbb F, s) +\delta_{[Z(s)]} =\mathsf{e}(\mathbb F, s)\,,
    $$
    where $Z(s) = s^{-1}(0)$. In particular, if $s$ has no zeroes on $M$, then
    $$
    \ddc\mathfrak g(\mathbb F, s) = \mathsf{e}(\mathbb F, s)
    $$
\end{lemma}
\begin{proof}
    By Lemma \ref{lem:ex-green}, this follows from \cite[Proposition 2.1.7 (3)]{GarciaSankaran2019}.
\end{proof}

\subsection{Moduli spaces of abelian varieties}

Let $\Lambda_g = \ZZ^{2g}$, which carries a canonical symplectic form given by the matrix
$$
\omega_g = \begin{pmatrix}
    0 & I_g\\
    -I_g & 0
\end{pmatrix}\,.
$$
The Siegel upper half-space is the symmetric domain for the symplectic group $\Sp_{2g}$:
$$
\mathfrak H_g = \{Z \in \operatorname{Mat}_{g\times g}(\mathbb C)\, : Z = Z^t\,, \ \Im(Z) >0\} = \Sp_{2g}(\RR)/ U(g)\,.
$$
We can also see the Siegel upper half-spaces as moduli spaces of marked abelian varieties:
$$
\mathfrak H_g = \{(A,\theta, \mu)\,: (A,\theta)\text{ is a principally polarized abelian variety, }\mu : (\Lambda_g, \omega_g)  \stackrel{\sim}\to (H_1(A, \ZZ), c_1(\theta))\}\,.
$$
Here, $Z$ corresponds to $\CC^g / \langle Z \mid 1\rangle$. See \cite[Section 6]{BL04}. The group of symplectic automorphisms of $(\Lambda_g, \omega_g)$ is $\operatorname{Sp}_{2g}(\ZZ)$, so it acts naturally on $\mathfrak H_g$, and we have a natural isomorphism of orbifolds\footnote{The reader who is not familiar with the notion of orbifold can take a level cover $\mathcal A_g(N) \to \mathcal A_g$, $N \geq 3$, perform calculations on the complex manifold $\mathcal A_g(N)$ and descend the consequences to $\mathcal A_g$. For a reference on differential forms, currents, vector bundles on orbifolds, see \cite{Behrend2004, SolomonTukachinsky2024}.}
$$
\mathcal A_g \cong [\mathfrak H_g / \Sp_{2g}(\ZZ)]\,.
$$
The \emph{Hodge bundle} is the trivial vector bundle
$$
\mathbb E_g  = \mathcal O_{\mathfrak H_g}^{\oplus g} = \mathbb C^g \times \mathfrak H_g\,.
$$
We equip it with the natural Hermitian metric given by the symmetric matrix $\Im (Z)$. 

We consider the complex manifold $M = \mathfrak H_1 \times \mathfrak H_g$ with coordinates
\begin{equation}\label{eq:coordinates}
    z = x+iy\,, \quad Z=X+iY\,.
\end{equation}

The trivial holomorphic bundle $ \mathbb F :=\mathbb E_1^\vee \otimes \mathbb E_{g}^\vee =\mathcal O_{\mathfrak H_1} \otimes \mathcal O_{\mathfrak H_g}^{\oplus g}$ is equipped with the metric given by the symmetric matrix $H^{-1}=(2yY)^{-1}$. Therefore, its connection is given by
\begin{equation}\label{eq:connection}
\nabla_{\mathbb F} = d +H\partial H^{-1} = d- \frac{dz}{2iy} - \frac{1}{2i}dZY^{-1}\, .
\end{equation}
Its curvature matrix is
\begin{equation}\label{eq:curvature}
    R:=\nabla_{\mathbb F}^2= \overline{\partial}(H\partial H^{-1})= -\frac{dz\wedge d\overline{z}}{4y^2}\operatorname{id}_g -\frac{1}{4}dZ Y^{-1}\wedge d\overline{Z}Y^{-1} \, .
\end{equation}
Similarly,
$$
\nabla_{\mathbb F^\vee} = d +H^{-1}\partial H = d+\frac{dz}{2iy} + \frac{1}{2i}Y^{-1}dZ
$$

We consider the variation of Hodge structures of weight $0$ given by $\mathbb L = \Hom (\Lambda_1 , \Lambda_g)$. We see elements $\ell \in \mathbb L$ as matrices
$$
\ell  = \begin{pmatrix}
    P_1& P_2\\Q_1&Q_2
\end{pmatrix}\in M_{2g\times 2}(\mathbb Z)
$$
The lattice $\mathbb L$ comes with a symmetric bilinear form
$$
\mathsf B (\ell_1, \ell_2) = \operatorname{tr}(\omega_1^{-1}\ell_1^t\omega_g \ell_2)\,,
$$
induced from $\omega_1$ and $ \omega_g$. Let $\mathsf{Q}(\ell) = \mathsf{B}(\ell, \ell)$ be the associated quadratic form. If we write $\ell = \begin{pmatrix}
    P_1& P_2\\Q_1&Q_2
\end{pmatrix}$, then
\begin{equation}\label{eq:norm-explicit}
\ell^t \omega_g \ell = \begin{pmatrix}
    0 & P_1^tQ_2 - Q_1^t P_2\\
    Q_1^t P_2 - P_1^tQ_2 & 0
\end{pmatrix} = (P_1^tQ_2 - Q_1^tP_2) \omega_1\,,
\end{equation}
so $\mathsf{Q}(\ell) = 2 P_1^t Q_2 - 2Q_1^t P_2$, and $(\mathbb L, \mathsf{Q})$ is isomorphic to the direct sum of $2g$ hyperbolic planes.

\begin{definition}
For any $\ell\in \Hom(\Lambda_1, \Lambda_g)$, we consider the section
$$
s_{\ell} : M \longrightarrow \mathbb F := \mathbb E^\vee_1 \otimes \mathbb E^\vee_g\, ,\quad s_{\ell}(z,Z) = z(ZP_2 + Q_2) -(ZP_1 + Q_1)\,.
$$
\end{definition}
Note that the definition also makes sense for $\ell \in \Hom(\Lambda_1, \Lambda_g)_{\RR}$.
\begin{lemma}\label{lem:section-equivariance}
    Let
    $$
    \gamma = (\gamma_1, \gamma_g) = \left(\begin{pmatrix}a&b\\c&d\end{pmatrix} ,\begin{pmatrix}A&B\\C&D\end{pmatrix}\right)  \in\SL_2(\RR) \times \Sp_{2g}(\RR)
    $$
    and $\ell\in\mathbb L$. Then, we have\footnote{The factor of automorphy $(cz + d)^{-1}(CZ+D)^{-t}$ is also the one used to descend $\mathbb E_1^\vee \otimes \mathbb E_g^\vee$ to a vector bundle on $\Ac_1 \times \Ac_g$.}
    $$
    s_{\gamma_g^{-t}\ell \gamma_1^t}(\gamma \cdot (z, Z)) = (cz+d)^{-1}(CZ+D)^{-t}s_{\ell}(z, Z)\,, \quad \mathsf{Q}(\gamma_g^{-t}\ell \gamma_1^t) = \mathsf{Q}(\ell)\,, \quad  \gamma^*\mathsf{e}_{\lambda}(\mathbb F, s_{\gamma_g^{-t}\ell \gamma_1^t}) = \mathsf{e}_\lambda(\mathbb F, s_\ell)\,.
    $$
\end{lemma}
\begin{proof}
    For the second identity, we have
    $$
    \begin{aligned}
    \mathsf{Q}(\gamma_g^{-t}\ell \gamma_1^t) &= \operatorname{tr} (\omega_1^{-1}\gamma_1\ell^t\gamma_g^{-1} \omega_g\gamma_g^{-t}\ell \gamma_1^t) & \text{(by definition of $\mathsf{B}$)}\\
    &=\operatorname{tr} (\omega_1^{-1}\gamma_1\ell^t\omega_g\ell \gamma_1^t)& \text{(since $\gamma_g \in \Sp_{2g}(\RR)$)}\\
    &=\operatorname{tr} (\gamma_1^t\omega_1^{-1}\gamma_1\ell^t\omega_g\ell )& \text{(since $\operatorname{tr}(AB) = \operatorname{tr}(BA)$)}\\
    &= \operatorname{tr} (\omega_1^{-1}\ell^t\omega_g\ell ) & \text{(since $\gamma_1 \in \SL_{2}(\RR)$)}
    \end{aligned}
    $$
    The holomorphic transformation $(cz + d)^{-1}(CZ+D)^{-t}$ preserves the Hermitian metric $H^{-1} = (2yY)^{-1}$ of $\mathbb F$, so it also preserves the Chern connection, the curvature and the norm of $s_\ell$. In particular 
    $$
    (\gamma_*s_\ell)(\gamma\cdot(z,Z)) = (cz+d)^{-1}(CZ+D)^{-t}s_\ell(z,Z)\,, \quad \gamma^*\mathsf{e}_{\lambda}(\mathbb F,\gamma_*s_\ell) = \mathsf{e}_\lambda(\mathbb F,s_\ell)\,,
    $$
    which gives the third identity after applying the first one.
\end{proof}
\begin{lemma}\label{lem:reg-section}
    If $\ell \neq 0$ then $s_{\ell}$ is a regular section, and its zero locus $Z(s_{\ell})$ parametrizes triples $(A, \theta_A, \mu_A, E, \theta_E, \mu_E, f)$, where $(A,\theta_A, \mu_A)$ is a marked principally polarized abelian variety of dimension $g$, $(E, \theta_E, \mu_E)$ is an elliptic curve, and $f :E \to A$ is a homomorphism such that
    $$
    \ell = \mu_A^{-1} \circ f_* \circ \mu_E : \Lambda_1 \to \Lambda_g\,.
    $$
    If $\mathsf Q(\ell) \leq 0$ then $Z(s_{\ell})$ is empty. If $\mathsf Q(\ell) > 0$ then the map $Z(s_{\ell}) \to \mathfrak H_g$ is proper.
\end{lemma}
\begin{proof}
To see why $s_\ell$ is a regular section, take its derivative:
$$Ds_{\ell}|_{(z,Z)}(u,U) = u(ZP_2 + Q_2) + U(zP_2-P_1).$$
Assume that $s_\ell(z,Z)=0$ and $\ell\neq 0$. We want to show that $Ds_\ell|_{(z,Z)}$ is a surjective linear map. Let $w = zP_2 - P_1\in \bbC^g$. If $w=0$, then since $z\in \mathfrak H$ has positive imaginary part, $P_1=P_2=0$. Using $s_\ell(z,Z)=0$, we get that $Q_1=Q_2=0$ also, which contradicts $\ell\neq 0$. Therefore, $w\neq 0$. Now, observe that $Ds_\ell|_{(z,Z)}(0,U) = Uw$ is a surjective linear map $\Sym^2(\bbC^g) \to \bbC^g$.

The zero locus $Z(s_\ell)$ parametrizes morphisms $f:E \to A$ of marked principally polarized abelian varieties; the induced map on homology is given by $\ell$, by Lemma 4.3 of \cite{GreerLian}, where $\ell$ is called $B_h$ instead. Observe that $\mathsf{Q}(\ell) = \operatorname{tr}(\theta_E\circ f^\vee\circ \theta_A \circ f|H_1(E,\bbZ))$, which is always a positive even integer. 

We now check that $Z(s_\ell) \to \mathfrak H_g$ is a closed holomorphic immersion, so in particular it is proper. Let $Z\in \mathfrak H_g$ be in the image of $Z(s_\ell)$, so there exists some $z\in \mathfrak H$ such that 
$$ZP_1+Q_1 = z(ZP_2+ Q_2).$$
This $z$ would be unique, provided that $ZP_2+Q_2\neq 0$. Assuming otherwise, we would have that $0=ZP_1+Q_1=ZP_2+Q_2$. Since $Z$ has positive definite imaginary part, this is only possible if $P_1=Q_1=P_2=Q_2=0$, which contradicts $\ell\neq 0$.
\end{proof}

For later use, we record the covariant derivatives of $s_\ell$.
\begin{align}
\begin{split}\label{eq:covariant-1}
    \nabla_{\mathbb F}^{1,0}s_{\ell}=\, &  \partial s_{\ell} + H \partial H^{-1}s_{\ell}\\
    =\,&\Big((ZP_2 +Q_2) - \frac{s_{\ell}}{2iy}\Big)dz + dZ(zP_2 -P_1) -\frac{1}{2i}dZ\,Y^{-1}s_{\ell}\\
    =\,&-\frac{1}{2i}\Big(\frac{s_{\ell}(\overline z,Z)}{y}\,dz + dZ\,Y^{-1}s_{\ell}(z,\overline Z)\Big)
\end{split}\\
\begin{split}\label{eq:covariant-2}
    \nabla_{\mathbb F^\vee}^{0,1} s_{\ell}^\dagger =\, &\overline{\nabla^{1,0}_{\mathbb F}s_{\ell}}^{\,t}\,(2yY)^{-1}\\
    =\,&\frac{1}{4iy}\Big(\frac{s_{\ell}(z,\overline Z)^t}{y}\,d\overline z + s_{\ell}(\overline z,Z)^t\,Y^{-1}d\overline Z\Big)Y^{-1}
\end{split}
\end{align}

% \begin{lemma}\label{lem:norm-and-covariant}
%     The norm and covariant derivatives of $s_{\ell}$ are given by 
%     {
%     \allowdisplaybreaks
%     \begin{align}
%         \begin{split}
%             \|s_{\ell}\|^2 =\,& \overline{s_{\ell}}^{\,t}(2yY)^{-1}s_{\ell} \\
%             =\,&\frac{1}{2y}\bigg(|z|^2(\overline{Z}P_2 + Q_2)^tY^{-1}(ZP_2 + Q_2) - 2 \operatorname{Re}\big(\overline{z}(\overline{Z}P_2 +Q_2)^tY^{-1}(ZP_1 + Q_1)\big) \\&+ (\overline{Z}P_1 + Q_1)^tY^{-1}(ZP_1 + Q_1)\bigg)
%         \end{split}\\
%         \begin{split}
%             \nabla_{\mathbb F}^{1,0}s_{\ell}=\, &  \partial s_{\ell} + H \partial H^{-1}s_{\ell}\\
%             =\,&\Big((ZP_2 +Q_2) - \frac{s_{\ell}}{2iy}\Big)dz + dZ(zP_2 -P_1) -\frac{1}{2i}dZ\,Y^{-1}s_{\ell}\\
%             =\,&-\frac{1}{2i}\Big(\frac{s_{\ell}(\overline z,Z)}{y}\,dz + dZ\,Y^{-1}s_{\ell}(z,\overline Z)\Big)
%         \end{split}\\
%         \begin{split}
%             \nabla_{\mathbb F^\vee}^{0,1} s_{\ell}^\dagger =\, &\overline{\nabla^{1,0}_{\mathbb F}s_{\ell}}^{\,t}\,(2yY)^{-1}\\
%             =\,& \left(d\overline{Z}\left(\overline{z}P_2 - P_1\right)+\frac{1}{2i}d\overline{Z}\,Y^{-1}\overline{s_{\ell}} +\Big(\frac{1}{2iy}\overline{s_{\ell}}+\overline{Z}P_2+Q_2\Big)d\overline{z}\right)^t\frac{1}{2y}Y^{-1}\\
%             =\,&\frac{1}{4iy}\Big(\frac{s_{\ell}(z,\overline Z)^t}{y}\,d\overline z + s_{\ell}(\overline z,Z)^t\,Y^{-1}d\overline Z\Big)Y^{-1}
%         \end{split}
%     \end{align}}
% \end{lemma}

\subsection{The Kudla--Millson form}

As a function of the modular parameter $\tau = u + i v$, the Kudla--Millson form is
\begin{equation}\label{eq:define-KM}
    \Theta_{\KM}(\tau) = \sum_{\ell \in \Hom(\Lambda_1, \Lambda_g)} e^{\pi i\tau \mathsf Q(\ell)} \mathsf{e}(\mathbb F, s_{v^{1/2}\ell})
\end{equation}
By \eqref{eq:norm-explicit}, the lattice $(\Hom(\Lambda_1, \Lambda_g), \mathsf Q)$ is even and unimodular of signature $(2g,2g)$ so we have
\begin{theorem}[Kudla--Millson]\label{thm:original-KM}
    If we see each summand in \eqref{eq:define-KM} as a smooth map $\mathscr{H} \times M \to \Omega^{g,g}(M)$, then \eqref{eq:define-KM} converges, uniformly on compact subsets of $\mathscr{H} \times M$. As a function of $(z,Z) \in M$, it is closed and $\Sp_{2g}(\ZZ) \times \Sp_2(\ZZ)$-invariant. As a function of $\tau \in \mathscr{H}$, it satisfies the modular equations
    $$
    \Theta_{\KM}(\tau + 1) = \Theta_{\KM}(\tau)\,, \quad \Theta_{\KM}(-1/\tau) = \tau^{2g}\Theta_{\KM}(\tau)\,.
    $$
\end{theorem}

The computation of the Fourier coefficients of the cohomology class $[\Theta_{\KM}(\tau)]$ is also part of the Kudla--Millson theorem. However, we are interested in the refinement that keeps track of the transgression relating them. This is the main advantage of working with the Mathai--Quillen formalism.
\begin{theorem}\label{thm:original-KM-green}
    The sum of currents (cf. Lemmas \ref{lem:green-is-current} and \ref{lem:reg-section})
    $$
    \mathfrak G (\tau) = \sum_{\ell \in \Hom(\Lambda_1, \Lambda_g) \setminus \{0\}} e^{\pi i \tau \mathsf Q(\ell)} \mathfrak g(\mathbb F, s_{v^{1/2}\ell})
    $$
    converges uniformly on compact subsets of $\mathscr{H} \times M$ and defines an $\SL_2(\ZZ)\times \Sp_{2g}(\ZZ)$-invariant current. Moreover, we have
    \begin{equation}\label{eq:to-be-extended}
            \ddc \mathfrak G(\tau) + \sum_{\ell \in \Hom(\Lambda_1, \Lambda_g) \setminus \{0\}} e^{\pi i \tau \mathsf{Q}(\ell)}\delta_{[Z(s_{\ell})]} = \Theta_{\KM}(\tau) - \mathsf{e}(\mathbb F, 0)
    \end{equation}
    as currents on $\mathscr{H} \times M$.
\end{theorem}

\subsection{The space $W_2$} The weight 2 polarized variation of Hodge structure $\mathbb L = \Hom(\Lambda_1,\Lambda_g)(-1) \simeq \Lambda_1^* \otimes \Lambda_g$ over $\mathcal A_1\times \mathcal A_g$ degenerates along the smooth divisor $i\infty \times \mathcal A_g$. This degeneration gives rise to a filtration of $\mathbb{L}$ called the weight filtration. To define it, let $$\Gamma_\infty  = \left\lbrace \left. \begin{pmatrix}
    1 & n\\ 0&1
\end{pmatrix} \right| n \in \mathbb{Z} \right\rbrace\subset \mathrm{PSL}_2(\ZZ)
$$ 
be the subgroup that stabilizes the cusp $i\infty$. Its generator $T=\left(\begin{smallmatrix}
1 & 1 \\ 0 & 1
\end{smallmatrix}\right)$ is unipotent, with logarithm $N=\log T=\left(\begin{smallmatrix}
0 & 1 \\ 0 & 0
\end{smallmatrix}\right)$. The weight filtration on $\mathbb{L}$ is $0 \subset W_2 \subset \mathbb{L}$, where $W_2$ is the space of matrices 
$$
W_2 = \mathrm{Im}(N) \otimes \Lambda_g = \left\{ \left. \ell = \begin{pmatrix}
    P_1 & 0\\
    Q_1 & 0
\end{pmatrix} \right| P_1,Q_1 \in M_g(\mathbb{Z}) \right\}.
$$

For any $Z \in \mathfrak H_g$, let $\mathfrak q_{Z} : W_2 \to \RR$ be the function defined by
$$
\mathfrak q_{Z} : (P_1, Q_1) \mapsto \frac{1}{2}(\overline{Z}P_1 +Q_1)^tY^{-1}(ZP_1 +Q_1)  = \frac{1}{2}\big(Q_1^tY^{-1}Q_1 + P_1^t(Y+XY^{-1}X)P_1 + 2P_1^tXY^{-1}Q_1\big)\,.
$$
Its matrix representation is given by
\begin{equation}\label{eq:factorize-quadratic}
\frac{1}{2}\begin{pmatrix}
    Y + XY^{-1}X&XY^{-1}\\
    Y^{-1}X&Y^{-1}
\end{pmatrix}
=
\frac{1}{2}\begin{pmatrix}
    I & X\\0&I
\end{pmatrix}
\begin{pmatrix}
    Y&0\\0&Y^{-1}
\end{pmatrix}
\begin{pmatrix}
    I&0\\X&I
\end{pmatrix}\,,
\end{equation}
so $\mathfrak q_Z$ is a positive definite quadratic form on $W_2$.

\begin{lemma}\label{lem:simplified-norm}
    If $\ell = \begin{pmatrix}
        P_1 & P_2 \\Q_1 & Q_2
    \end{pmatrix}$ then we have
    $$
    \|s_\ell\|^2 + \frac12 \mathsf{Q}(\ell) = y \mathfrak q_Z(P_2, Q_2) + y^{-1}\mathfrak q_Z(P_1-xP_2, Q_1-xQ_2)
    $$
\end{lemma}
\begin{proof}
    If we write $s_\ell = iy(ZP_2 + Q_2) - (ZP_1 + Q_1 - x(ZP_2 + Q_2))$ and expand $\|s_\ell^2\| = \overline{s_\ell}^t(2yY)^{-1}s_\ell$, we obtain
    \begin{align*}
        \frac{1}{2y}\Big(&y^2 (\overline{Z}P_2 + Q_2)^{t}Y^{-1}(ZP_2 + Q_2) + (\overline{Z}P_1 + Q_1 - x(\overline{Z}P_2 + Q_2))Y^{-1}(ZP_1 + Q_1 - x(ZP_2 + Q_2))\\
        & \quad + iy\big((\overline{Z}P_2 + Q_2)^{t}Y^{-1}(ZP_1 + Q_1 - x(ZP_2 + Q_2)) - (\overline{Z}P_1 + Q_1 - x(\overline{Z}P_2 + Q_2))Y^{-1}(ZP_2 + Q_2)\big) \Big)\\
        =\,&  y \mathfrak q_Z(P_2, Q_2) + y^{-1}\mathfrak q_Z(P_1-xP_2, Q_1-xQ_2)\\
        &\quad + \frac{i}{2}\big((\overline{Z}P_2 + Q_2)^{t}Y^{-1}(ZP_1 + Q_1) - (\overline{Z}P_1 + Q_1)Y^{-1}(ZP_2 + Q_2)\big)\\
        =\,&  y \mathfrak q_Z(P_2, Q_2) + y^{-1}\mathfrak q_Z(P_1-xP_2, Q_1-xQ_2) - \operatorname{Im} \big((\overline{Z}P_2 + Q_2)Y^{-1}(ZP_1 + Q_1)\big)
    \end{align*}
    An explicit calculation shows that the imaginary part of $(\overline{Z}P_2 + Q_2)Y^{-1}(ZP_1 + Q_1)$ is $P_1^t Q_2 - Q_1^tP_2$, which is precisely $\frac{1}{2}\mathsf{Q}(\ell)$ by \eqref{eq:norm-explicit}.
\end{proof}
In particular, the Schmid estimate
$$
\|s_\ell\|^2 + \frac12 \mathsf{Q}(\ell) \geq Cy\,, \qquad \ell \not \in W_2
$$
clearly holds, where $C$ is uniform on compact subsets of $M$, and $\|s_\ell\|^2 = y^{-1}\mathfrak q_Z(\ell)$ for all $\ell \in W_2$.

\begin{lemma}\label{lem:homogeneous-euler-W2}
    For $\ell \in W_2$, $s_\ell$ is independent of $z=x+iy \in \mathfrak H_1$, and we can write
    \begin{align*}
    \mathsf{e}_{\lambda} (\mathbb F, s_\ell) =\,& \sum_{i,j\in \{0,1\}}e^{-2\pi\mathfrak{q}_Z(y^{-1/2}\ell)}E^{(i,j)}_{\lambda} ( y^{-1/2}\ell)\left(\frac{dz}{y}\right)^i\wedge \left(\frac{d\overline{z}}{y}\right)^j\\
    \mathrm{d}^c \mathsf{e}_{\lambda} ( \mathbb F, s_\ell) =\,& \sum_{i,j\in \{0,1\}}e^{-2\pi\mathfrak{q}_Z(y^{-1/2}\ell)}F^{(i,j)}_{\lambda} ( y^{-1/2}\ell)\left(\frac{dz}{y}\right)^i\wedge \left(\frac{d\overline{z}}{y}\right)^j\,,
    \end{align*}
    where $E^{(i,j)}_{\lambda}, F^{(i,j)}_\lambda : W_2 \to A^{*-i,*-j}(\mathfrak H_g)[\lambda]$ are polynomials.
\end{lemma}
\begin{proof}
    For $\ell \in W_2$, the section $s_\ell(z,Z)$ only depends on $Z$. By \eqref{eq:covariant-1} and \eqref{eq:covariant-2}, we can write
    \begin{align*}
    \nabla^{1,0}s_\ell \wedge \nabla^{0,1}s_\ell^\dagger =\, & \frac{1}{8y}\left(s_\ell (Z)\frac{dz}{y} \wedge s_{\ell}(\overline{Z})^t\frac{d\overline{z}}{y} + dZY^{-1}s_{\ell}(\overline{Z})\wedge s_\ell(\overline{Z})^t \frac{d\overline{z}}{y}\right.\\
    & \left.\qquad  + s_\ell (Z)\frac{dz}{y}\wedge s_\ell(Z)^tY^{-1}d\overline{Z} + dZY^{-1}s_{\ell}(\overline{Z})\wedge s_\ell(Z)^tY^{-1}d\overline{Z}\right) Y^{-1}\\
    =\,& \sum_{i,j\in \{0,1\}} B^{(i,j)}(y^{-1/2}\ell) \left(\frac{dz}{y}\right)^i\wedge \left(\frac{d\overline{z}}{y}\right)^j\,,
    \end{align*}
    where $B^{(i,j)}(\xi)$ is polynomial in $\xi$. The curvature can also be written as
    $$
    R =\,-\frac{1}{4}\frac{dz}{y}\wedge \frac{d\overline{z}}{y}\operatorname{id}_g -\frac{1}{4}dZ Y^{-1}\wedge d\overline{Z}Y^{-1} = R^{(1,1)} \frac{dz}{y}\wedge \frac{d\overline{z}}{y} + R^{(0,0)}
    $$
The first formula in the statement of the lemma follows by applying the definition of $\mathsf{e}_\lambda$.

The operator $\mathrm{d}^c$ preserves this structure, since
$$
\mathrm{d}^c\left(\frac{dz}{y}\right) = -\frac{1}{8\pi}\,\frac{dz}{y}\wedge\frac{d\overline z}{y}\,, \quad \mathrm{d}^c\left(\frac{d\overline{z}}{y}\right) = \frac{1}{8\pi}\,\frac{dz}{y}\wedge\frac{d\overline z}{y}\,,
$$
and, if we abbreviate $\xi = y^{-1/2}\ell$, then
\begin{equation*}
\mathrm{d}^c\Big(e^{-2\pi \mathfrak q_Z(\xi)}E_\lambda^{(i,j)}(\xi)\Big)
= \mathrm{d}^c_Z\Big(e^{-2\pi \mathfrak q_Z(\xi)}E_\lambda^{(i,j)}(\xi)\Big)+ \frac{1}{16\pi}\Big(\frac{dz}{y}+\frac{d\overline z}{y}\Big)\wedge \Big(\sum_{k=1}^{2g} \xi_k\frac{\partial}{\partial \xi_k}\Big)\Big(e^{-2\pi \mathfrak q_Z(\xi)}E_\lambda^{(i,j)}(\xi)\Big)\,.
\end{equation*}
\end{proof}

The next proposition is the key evaluation needed for the proof of Theorem \ref{thm:main} in the next section.

\begin{proposition}\label{p:W2exact0}
    Let $\lambda$ be a formal parameter, and let $\mathrm{d}\ell$ be the Lebesgue measure on $W_{2, \RR}$. Then
    \begin{equation}\label{eq:the-integral}
    \int_{W_{2, \mathbb R}}\mathsf{e}_\lambda (\mathbb F, s_{\ell})\, d\ell = (y\lambda)^g\, .
    \end{equation}
\end{proposition}
\begin{proof}
    Let $P_{\infty} = \Gamma_{\infty} \times \Sp_{2g}$ be the parabolic that stabilizes the cusp $i \infty$, and let
    $$
    \gamma = (\gamma_1, \gamma_g) = \left(\begin{pmatrix}
        a&b\\
        0&a^{-1}
    \end{pmatrix}, \begin{pmatrix}
        A&B\\C&D
    \end{pmatrix}\right)\in \SL_2(\RR) \times \Sp_{2g}(\RR)\,,
    $$
    be an element of $P_\infty(\RR)$. Then, by Lemma \ref{lem:section-equivariance} and the change of variables, we have
    \begin{align*}
    \gamma^* \int_{W_{2, \mathbb R}}\mathsf{e}_\lambda (\mathbb F, s_{\ell})\, d\ell =\,& \int_{W_{2, \mathbb R}}\gamma^*\mathsf{e}_\lambda (\mathbb F, s_{\ell})\, d\ell\\
    =\, & \int_{W_{2,\RR}}\mathsf{e}_{\lambda} (\mathbb F, s_{\gamma^{-1}\cdot \ell}) d\ell\\
    =\, & \det (\gamma : W_2 \to W_2)\int_{W_{2,\RR}}\mathsf{e}_{\lambda} (\mathbb F, s_{\ell}) d\ell\\
    =\,& a^{2g} \int_{W_{2,\RR}}\mathsf{e}_\lambda(\mathbb F, s_\ell)d\ell\,.
    \end{align*}
    Therefore, both sides of \eqref{eq:the-integral} transform accordingly under the action of $P_{\infty}(\RR)$. 
    
    Since the action of $P_\infty(\RR)$ on $M$ is transitive, it is enough to evaluate \eqref{eq:the-integral} at $(z,Z) = (i, i\operatorname{id}_g)$. Put $s_i=q_i+ip_i$, so that $s_{\ell}=-iP_1-Q_1$. Then the norm of $s_\ell$, its covariant derivatives \eqref{eq:covariant-1}, \eqref{eq:covariant-2} and the curvature become
    \begin{align}
    \begin{split}\label{eq:at-ii-i}
    \lVert s_{\ell}\rVert^2=&\frac{1}{2}\sum_{i=1}^g|s_i|^2\,,\\
     (\nabla^{1,0}s_{\ell})_i=&\,\frac{i}{2}s_i\,dz+\frac{i}{2}\sum_k dZ_{ik}\,\overline{s}_k\,,\\
     (\nabla^{0,1}s_{\ell}^\dagger)^j
     =&\,-\frac{i}{4}\overline{s}_j\,d\overline z
       -\frac{i}{4}\sum_k d\overline Z_{jk}\,s_k\,,\\
     R_{ij}
     =&\,-\frac{\delta_{ij}}{4}\,dz\wedge d\overline z
       -\frac14\sum_k dZ_{ik}\wedge d\overline Z_{kj}\,.
       \end{split}
    \end{align}
    If we expand the definition of $\mathsf{e}_\lambda$, with notation $R$ and $B_s$ as in Lemma \ref{lem:ex-euler}, we find that
    $$
    \mathsf{e}_\lambda(\mathbb F, s_\ell) = e^{-\pi\sum_i |s_i|^2}\sum_{K\subset\{1,\ldots,g\}}\lambda^{g-|K|}\Big(\frac{-1}{2\pi i}\Big)^{|K|}
    \sum_{I\subset K}\frac{(2\pi)^{|I|}}{|I|!}
    \sum_{\sigma\in\mathfrak S_K}\operatorname{sgn}(\sigma)
    \prod_{i\in K\smallsetminus I}R_{i\sigma(i)}\prod_{i\in I}(B_s)_{i\sigma(i)}\,.
    $$
    The curvature terms do not depend on $\ell \in W_2$. To evaluate the $B_s$ terms, we introduce the Gaussian average
    $$
    \langle f\rangle=\int_{W_{2,\mathbb R}}e^{-\pi\sum_i|s_i|^2}f\,d\ell.
    $$
    Then, using the explicit formula in \eqref{eq:at-ii-i}, one sees that
    $$
    \langle 1\rangle=1\,,\qquad \langle s_i\overline s_j\rangle=\frac{1}{\pi}\delta_{ij}\,, \qquad \langle s_is_j\rangle= \langle\overline s_i\overline s_j\rangle=0\, .
    $$
    Consequently,
    $$
     \left\langle(\nabla^{1,0}s_{\ell})_i(\nabla^{1,0}s_{\ell})_j\right\rangle=0\,,\qquad\left\langle(\nabla^{0,1}s_{\ell}^\dagger)^i(\nabla^{0,1}s_{\ell}^\dagger)^j\right\rangle=0\, ,\qquad\left\langle(\nabla^{1,0}s_{\ell})_i(\nabla^{0,1}s_{\ell}^\dagger)^j\right\rangle=-\frac{R_{ij}}{2\pi}\, .
    $$
    Using Wick's theorem for $1$-forms $x_1, \ldots, x_{2r}$ whose coefficients are linear in $\ell$:
    $$
    \langle f_1\cdots f_{2r} \rangle = \sum_{\{i_1, j_1, \ldots, i_r, j_r\}\to [2r]}\operatorname{sgn}(\eps)\,\langle f_{i_1}f_{j_1}\rangle\cdots \langle f_{i_r}f_{j_r}\rangle\,,
    $$
    we have
    $$
    \left\langle\prod_{i\in I}(B_s)_{i \sigma(i)}\right\rangle = \left\langle\prod_{i\in I}(\nabla^{1,0}s_\ell)_i(\nabla^{0,1}s^\dagger_\ell)^{\sigma(i)}\right\rangle = \sum_{\tau : I \to \sigma(I)}\Big(\frac{-1}{2\pi}\Big)^{|I|}(-1)^{\sigma|_I^{-1}\circ \tau}\prod_{j\in I} R_{j\tau(j)}\,,
    $$
    where the sum runs over bijections $\tau$. Therefore,
    $$
    \int_{W_{2,\mathbb R}}\mathsf{e}_\lambda(\mathbb F, s_\ell)\, d\ell = \sum_{K\subset\{1,\ldots,g\}}\lambda^{g-|K|}\Big(\frac{-1}{2\pi i}\Big)^{|K|} \sum_{\substack{\sigma \in \mathfrak S_K\\I \subset K}}\sum_{\tau :I \to \sigma(I)}\frac{(-1)^{|I|}(-1)^{\sigma}(-1)^{\sigma|_I^{-1}\circ \tau}}{|I|!}\prod_{i \in K\smallsetminus I}R_{i \sigma(i)}\prod_{i \in I}R_{i \tau(i)}\,.
    $$
    If $\rho \in \mathfrak S_K$ is $\sigma$ outside of $I$ and $\tau$ on $I$, then $(-1)^\rho = (-1)^{\sigma}(-1)^{\sigma|_I^{-1}\circ \tau}$, and we sum over the ways of writing $\rho$ in terms of a $\sigma$ and a $\tau$, which are determined by a bijection $\sigma_I : I \to \rho(I)$:
    \begin{align*}
    \int_{W_{2,\mathbb R}}\mathsf{e}_\lambda(\mathbb F, s_\ell)\, d\ell =&\, \sum_{K\subset\{1,\ldots,g\}}\lambda^{g-|K|}\Big(\frac{-1}{2\pi i}\Big)^{|K|}\sum_{\rho\in \mathfrak S_K}\sum_{I\subset K}\sum_{\sigma_I : I \to \rho(I)}\frac{(-1)^{|I|}}{|I|!}(-1)^\rho\prod_{i\in K}R_{i\rho(i)}\\
    =&\, \sum_{K\subset\{1,\ldots,g\}}\lambda^{g-|K|}\Big(\frac{-1}{2\pi i}\Big)^{|K|}\det(R_{KK})\sum_{I\subset K}(-1)^{|I|}\\
    =&\, \sum_{K\subset\{1,\ldots,g\}}\lambda^{g-|K|}\Big(\frac{-1}{2\pi i}\Big)^{|K|}\det(R_{KK})(1-1)^{|K|}\\
    =&\, \lambda^g\,.
    \end{align*}
\end{proof}
\section{Proof of Theorem \ref{thm:main}}\label{sec:proof}

The proof is divided into several steps. We follow the method in \cite{Garcia}; that is, we analyze the \emph{Green equation} \eqref{eq:to-be-extended} against test forms which are not necessarily compactly supported on $\Ac_1 \times \Ac_g$.
The first goal of this section will be to show the following:
\begin{theorem}\label{thm:integrability}
    The Kudla--Millson form $\Theta_{\KM}(\tau)$ extends naturally to a current on $\Abar_1 \times \Ac_g$, which we denote by $\overline{\Theta}_{\KM}(\tau)$. As a current, it satisfies the modular equations
    $$
    \overline{\Theta}_{\KM}(\tau+1) = \overline{\Theta}_{\KM}(\tau)\,, \quad \overline{\Theta}_{\KM}(-1/\tau) = \tau^{2g}\overline{\Theta}_{\KM}(\tau)\,. 
    $$
\end{theorem}

\begin{notation}\label{notation}
    We will use the local coordinate $p= e^{2\pi i z}$ near the cusp. Currents are linear functionals on the space of \emph{test forms}. In this section, we will use $\eta$ to denote test forms on $\Abar_1 \times \Ac_g$ or in the parabolic neighborhood $[\mathfrak H_1/\Gamma_\infty] \times \Ac_g$. The restriction of $\eta$ to the universal cover $M$ can be written as
    $$
    \eta = \eta^{(0,0)} + \eta^{(1,0)} dz + \eta^{(0,1)}d\overline{z} + \eta^{(1,1)}dz\wedge d\overline{z}\,.
    $$
    In terms of the cusp coordinate, we have
    $$
    \eta = \eta^{(0,0)} + \eta^{(1,0)} \frac{dp}{2\pi i p} - \eta^{(0,1)}\frac{d\overline{p}}{2\pi i \overline{p}} + \eta^{(1,1)}\frac{dp\wedge d\overline{p}}{4\pi^2 |p|^2}\,,
    $$
    so for $\eta$ to extend smoothly to the cusp, we need the bounds
    \begin{equation}\label{eq:eta-decay}
    |\eta^{(0,0)}| = O(1)\,, \quad |\eta^{(1,0)}| + |\eta^{(0,1)}| = O(e^{-\pi y})\,, \quad |\eta^{(1,1)}| = O(e^{-2\pi y})\,.
    \end{equation}
    (and the same bounds should also hold for all derivatives of $\eta$).
\end{notation}

For $\eps \leq \eps_0 :=e^{-2\pi}$, let
$$
\mathsf{K}_{i\infty}(\eps) = \Big\lbrace x+iy \in \mathfrak H_1: y > -\frac{\log \eps}{2\pi}\Big\rbrace\,.
$$
Note that $\mathsf{K}_{i\infty}(\eps)$ is $\Gamma_\infty$-invariant, and defines a punctured disk in $[\mathfrak H_1 / \Gamma_\infty]$. The induced map $[\mathsf{K}_{i\infty}(\eps)/\Gamma_\infty] \to \Abar_1$ is a biholomorphism onto the $\eps$-neighbourhood of the cusp in the coordinate $p$:
\begin{equation}\label{eq:disk-around-cusp}
[\mathsf{K}_{i\infty}(\eps) /\Gamma_\infty] \stackrel{\sim} \longrightarrow \mathsf{B}_{i\infty}(\eps) \smallsetminus \{i \infty\}\,, \quad \mathsf{B}_{i\infty}(\eps) :=\{p = e^{2\pi i z} :  |p|< \eps\}
\end{equation}
We can write any compactly supported form on $\Abar_1 \times \Ac_g$ as $\eta = \eta_1 + \eta_2$ where
\begin{equation}\label{eq:support}
\operatorname{supp}(\eta_1)\Subset \Ac_1 \times \Ac_g\,, \quad \operatorname{supp} (\eta_2) \Subset \mathsf{B}_{i\infty}(\eps_0) \times \Ac_g\,.
\end{equation}

Then, by Theorem \ref{thm:original-KM}, it is clear that $\Theta_{\KM}(\tau) \wedge \eta_1$ converges and satisfies the modularity in Theorem \ref{thm:integrability}. Therefore, it is enough to work with forms which are supported on $\mathsf{B}_{i\infty}(\eps_0)\times \mathcal A_g$.

\begin{lemma}\label{lem:convergence-deltas}
    The sum
    $$
    \sum_{\ell \in \Hom(\Lambda_1, \Lambda_g)\neq 0} e^{\pi i \tau \mathsf{Q}(\ell)}\delta_{Z(s_{\ell})} 
    $$
    converges to a  current on $\mathsf{B}_{i \infty}(\eps_0)\times \Ac_g$, uniformly over compact subsets of $\mathscr{H}$.
\end{lemma}
\begin{proof}
    Note that, by Lemma \ref{lem:reg-section}, it is enough to consider $\ell$ outside of $W_2$. If we differentiate $s_\ell =0$, we obtain
    $$
    (ZP_2 + Q_2)\d z = -\d Z(zP_2 - P_1) = y \d ZY^{-1}(\overline{Z}P_2 + Q_2)
    $$
    Since $\ell$ does not belong to $W_2$, we have $ZP_2 +Q_2 \neq 0$, so we obtain that $dz = O(y)$, uniformly in $\ell$, over compact subsets of $\mathfrak H_g$. By \eqref{eq:eta-decay}, we have
    \begin{equation}\label{eq:some-step}
    \int \delta_{Z(s)} \wedge \eta = \int_{Z(s)} \sum_{i,j \in \{0,1\}} \eta^{(i,j)}(dz)^i\wedge (d\overline{z})^j = O(1) + O(ye^{-\pi y}) + O(y^2e^{-2\pi y}) = O(1)\,.
    \end{equation}
    The integral is uniformly bounded by a constant $C(\eta)$ that does not depend on $\ell$, and is uniform over $\mathfrak H_g \times \mathfrak H_1$.

    By Lemma \ref{lem:simplified-norm}, the following holds on $Z(s_\ell)$:
    $$
    \mathsf{Q}(\ell) = 2y \mathfrak{q}_Z(P_2, Q_2) + \frac{2}{y}\mathfrak{q}_Z (P_1 - xP_2, Q_1 - xQ_2)
    $$
    The quadratic form $\mathfrak{q}_Z(P_2, Q_2) + \mathfrak{q}_Z (P_1 - xP_2, Q_1 - xQ_2)$ is positive definite in $\ell \in \Hom(\Lambda_1, \Lambda_g)$. Therefore, we find that
    $$
    \begin{aligned}
        \mathsf{Q}(\ell) \geq \,&   2y \mathfrak{q}_Z(P_2, Q_2) \geq cy\,,\\
        \mathsf{Q}(\ell) \geq \,& y^{-1} \mathfrak{q}_Z(P_2, Q_2) + \mathfrak{q}_Z (P_1 - xP_2, Q_1 - xQ_2) \geq c'y^{-1}|\ell|^2\,.
    \end{aligned}
    $$
    Here $c = \min_{(P_2, Q_2) \in \mathbb Z^{2g}\smallsetminus 0} \mathfrak q_Z(P_2, Q_2)$, $c'$ is some constant depending uniformly on $Z$, $x$ and $|\ell|^2$ is the Euclidean inner product on $\Hom(\Lambda_1, \Lambda_g) \cong \ZZ^{4g}$. Multiplying these two inequalities, we find that
    $$
    \mathsf{Q}(\ell) \geq c'' |\ell|\,,
    $$
    where $c''$ is uniform over $\mathfrak{H}_1 \times \mathfrak H_g$ and independent of $y$. In particular, by \eqref{eq:some-step}, we have
    $$
    \bigg|\int \sum_{\ell \in \Hom(\Lambda_1, \Lambda_g)\neq 0} e^{\pi i \tau \mathsf{Q}(\ell)}\delta_{Z(s_{\ell})} \wedge \eta \bigg| \leq \sum_{\ell \in \Hom(\Lambda_1, \Lambda_g)} e^{-\pi v c|\ell|} C(\eta) < \infty
    $$
\end{proof}

The space $W_2$ is $\Gamma_\infty$-invariant. Therefore, for any $\eta$ with $\operatorname{supp}(\eta)\subset \mathsf{B}_{i\infty}(\eps_0)\times \Ac_g$, we have
$$
\int \Theta_{\KM}(\tau)\wedge \eta = \int_{[\mathfrak H_1/\Gamma_\infty]\times \Ac_g} \sum_{\ell \in W_2} \mathsf{e}(\mathbb F,s_{v^{1/2}\ell})\wedge \eta + \int_{[\mathfrak H_1/\Gamma_\infty]\times \Ac_g}\sum_{\ell \in \Hom(\Lambda_1, \Lambda_g)\setminus W_2}e^{\pi i \tau \mathsf{Q}(\ell)}\mathsf{e} (\mathbb F, s_{v^{1/2}\ell})\wedge \eta\,,
$$
where \eqref{eq:disk-around-cusp} allows us to see $\eta$ as a form on $\Ac_1\times \Ac_g$ or $[\mathfrak H_1/\Gamma_\infty]\times \Ac_g$. Both sums are $\Gamma_\infty \times \Sp_{2g}(\ZZ)$-invariant by Lemma \ref{lem:section-equivariance}. A priori, we do not know if any of the terms here converges, but it is clear that if both terms on the right side converge, so does the left one, and the equality will hold.

% \begin{lemma}
%     For any $\ell \in \Hom(\Lambda_1, \Lambda_g)$, the Mathai--Quillen Euler class is integrable in $\tau = u+iv\in \mathscr{H}$, uniformly over compact subsets of $\mathscr{H}$:
%     $$
%     \int e^{\pi i \tau \mathsf{Q}(\ell)}\mathsf{e}(\mathbb F, s_{v^{1/2}\ell})\wedge \eta < \infty
%     $$
% \end{lemma}

% We write
% $$
% \Theta_{\KM}(\tau) = \sum_{\ell \in W_2} \mathsf{e}(\mathbb F, s_{v^{1/2}\ell}) + \sum_{\ell \in \Hom(\Lambda_1, \Lambda_g)\smallsetminus W_2} e^{\pi i \tau \mathsf{Q}(\ell)}\mathsf{e}(\mathbb F, s_{v^{1/2}\ell})
% $$
% and analyze each term separately.
Thus, we are reduced to proving the following two results
\begin{proposition}\label{prop:integrability-noW2}
    The sum
    $$
    \sum_{\ell \in \Hom(\Lambda_1, \Lambda_g)\smallsetminus W_2}  e^{\pi i \tau \mathsf{Q}(\ell)}\mathsf{e}(\mathbb F, s_{v^{1/2}\ell})
    $$
    converges to a current on $\mathsf{B}_{i\infty}(\eps_0)\times \mathcal A_g$ uniformly on compact subsets of $\mathscr{H}$.
\end{proposition}

\begin{proposition}\label{prop:integrability-W2}
    The sum
    $$
    \sum_{\ell \in  W_2}  e^{\pi i \tau \mathsf{Q}(\ell)}\mathsf{e}(\mathbb F, s_{v^{1/2}\ell})
    $$
    converges to a current on $\mathsf{B}_{i\infty}(\eps_0)\times \mathcal A_g$ uniformly on compact subsets of $\mathscr{H}$.
\end{proposition}

The proof is deferred to Sections \ref{s:integrability-outside-W2} and \ref{s:integrability-W2}. Once Propositions \ref{prop:integrability-noW2} and \ref{prop:integrability-W2} have been established, the modular equations from Theorem \ref{thm:integrability} follow, since the map from differential forms on $\Ac_1 \times \Ac_g$ to currents on $\Abar_1 \times \Ac_g$ is injective.

The second goal is to show that the Green equation \eqref{eq:to-be-extended} still holds, at least for $g \geq 2$:
\begin{theorem}\label{thm:new-green}
    There is a transgression current $\overline{\Psi}(\tau)$ on $\Abar_1 \times \Ac_g$ which extends the current $\Psi :=\mathrm{d}^c\mathfrak G(\tau)$ from Theorem \ref{thm:original-KM-green}. If $g\geq 2$, we have an equality of currents on $\Abar_1 \times \Ac_g$
    $$
    \overline{\Theta}_{\KM}(\tau) = \mathrm{d} \overline{\Psi}(\tau)+ 
        \mathsf{e}(\mathbb F,0) + \sum_{\ell \in \Hom(\Lambda_1, \Lambda_g)\setminus\{0\}}e^{\pi i \tau \mathsf{Q}(\ell)}\delta_{Z(s_\ell)}\,.
    $$
    If $g=1$, we have
    $$
    \overline{\Theta}_{\KM}(\tau) = \mathrm{d} \overline{\Psi}(\tau) + 
        \mathsf{e}(\mathbb F,0)  +  \frac{1}{4\pi v}\delta_{\{i \infty\}\times \Ac_1} + \sum_{\ell \in \Hom(\Lambda_1, \Lambda_1)\setminus\{0\}}e^{\pi i \tau \mathsf{Q}(\ell)}\delta_{Z(s_\ell)}\,.
    $$
    In particular,
    \begin{equation}\label{eq:power-series-inA1xAg}
    c_g(\mathbb E_{\Abar_1}^\vee \otimes \EE_{\Ac_g}^\vee)  + \sum_{d \geq 1}[\NL_{g,d}^{+}]q^d \in H^{2g}(\Abar_1 \times \Ac_g)\llbracket q\rrbracket
    \end{equation}
    is the holomorphic part of the Fourier expansion at $i \infty \in \mathscr{H}$ of a holomorphic (for $g \geq 2$, and almost holomorphic for $g=1$) modular form of weight $2g$ for $\SL_{2}(\ZZ)$.
\end{theorem}

Theorem \ref{thm:new-green} implies our main result:

\begin{proof}[Proof of Theorem \ref{thm:main}]
    If we integrate the holomorphic Fourier expansion \eqref{eq:power-series-inA1xAg} along the projection $\pi : \Abar_1 \times \Ac_g \to \Ac_g$, and use that
    $$
    \pi_* c_g(\mathbb E_{\Abar_1}^\vee \otimes \EE_{\Ac_g}^\vee) = (-1)^g\big(\pi_*(\lambda_1^{\Abar_1})\lambda_{g-1}+ \pi_*(1)\lambda_g\big) = \frac{(-1)^g}{24}\lambda_{g-1}\,, \quad \pi_*[\NL_{g,d}^+] = [\NL_{g,d}]\,,
    $$
    we obtain the statement of Theorem \ref{thm:main}.
\end{proof}

Similarly to how we did with Theorem \ref{thm:integrability}, we will split the test form $\eta$ as $\eta = \eta_1 + \eta_2$, where $\eta_1, \eta_2$ satisfy \eqref{eq:support}. Theorem \ref{thm:new-green}, applied to $\eta_1$, is covered by Theorem \ref{thm:original-KM-green}\footnote{Note that $c_1(\mathbb E_{\Abar_1})$ is exact on any compact subset of $\Ac_1$.}, so we can study instead forms supported on $\mathsf{B}_{i\infty}(\eps_0)$. Therefore, Theorem \ref{thm:new-green} will be a consequence of the following two propositions, one in $W_2$ and one outside of $W_2$:

\begin{proposition}\label{prop:exact-not-W2-green}
    The current
    $$
    \sum_{\ell \in \Hom(\Lambda_1, \Lambda_g)\smallsetminus W_2} e^{\pi i \tau\mathsf{Q}(\ell)}\mathrm{d}^c\mathfrak{g}(\mathbb F, s_{v^{1/2}\ell}) 
    $$
    converges to a current on $\mathsf{B}_{i\infty}(\eps_0)\times \mathcal A_g$, which we denote by $\overline{\Psi}_{W_2^c}(\tau)$. We have
    \begin{equation}\label{eq:Green-outside-W2}
        \sum_{\ell \in \Hom(\Lambda_1, \Lambda_g)\smallsetminus W_2}e^{\pi i \tau \mathsf{Q}(\ell)}\mathsf{e}(\mathbb F, s_{v^{1/2}\ell}) = \mathrm{d}\overline{\Psi}_{W_2^c}(\tau) + \sum_{\ell \in \Hom(\Lambda_1, \Lambda_g)\smallsetminus \{0\}} e^{\pi i \tau \mathsf{Q}(\ell)}\delta_{Z(s_\ell)}
    \end{equation}
    as currents on $[\mathsf{K}_{i\infty}(\eps_0)/\Gamma_\infty]\times \mathcal A_g$.
\end{proposition}

\begin{proposition}\label{prop:exact-W2-green}
    The current
    $$
    \sum_{\ell \in W_2 \smallsetminus \{0\}} \mathrm{d}^c\mathfrak{g}(\mathbb F, s_{v^{1/2}\ell}) 
    $$
    converges to a current on $\mathsf{B}_{i\infty}(\eps_0)\times \mathcal A_g$, which we denote by $\overline{\Psi}_{W_2}(\tau)$. If $g \geq 2$, we have
    \begin{equation*}
        \sum_{\ell \in W_2 \smallsetminus \{0\}}\mathsf{e}(\mathbb F, s_{v^{1/2}\ell}) = \mathrm{d}\overline{\Psi}_{W_2}(\tau)\,,
    \end{equation*}
    as currents on $\mathsf{B}_{i\infty}(\eps_0)\times \mathcal A_g$. If $g=1$, we have
    \begin{equation*}
        \sum_{\ell \in W_2 \smallsetminus \{0\}}\mathsf{e}(\mathbb F, s_{v^{1/2}\ell}) = \mathrm{d}\overline{\Psi}_{W_2}(\tau) + \frac{1}{4\pi v}\delta_{\{i \infty\} \times \Ac_1}\,,
    \end{equation*}
\end{proposition}

Their proofs are delayed to Sections \ref{s:gree-out-W2} and \ref{s:green-in-W2}. The transgression current $\overline{\Psi}(\tau)$ is obtained as
$$
\overline{\Psi}(\tau) = \overline{\Psi}_{W_2^c}(\tau) + \overline{\Psi}_{W_2}(\tau)\,.
$$
Combining Propositions \ref{prop:exact-not-W2-green} and \ref{prop:exact-W2-green}, we obtain the first part of Theorem \ref{thm:new-green}. By Theorem \ref{thm:integrability}, $\Theta_{\KM}$ satisfies the modular equations. Since $\mathrm{d}\overline{\Psi}(\tau)$ is exact, the Fourier expansion in $q= e^{2\pi i \tau}$ of the cohomology class $[\overline{\Theta}_{\KM}]$ is given by
$$
\mathsf{e}(\mathbb F,0) + \sum_{\ell \in \Hom(\Lambda_1, \Lambda_g)\setminus\{0\}}e^{\pi i \tau \mathsf{Q}(\ell)}\delta_{Z(s_\ell)} = \mathsf{e}(\mathbb F,0)  + \sum_{d\geq 1}q^{d}[\NL_{g,d}^+]
$$
if $g\geq 2$, which is holomorphic, and
$$
\mathsf{e}(\mathbb F,0)  +  \frac{1}{4\pi v}\delta_{\{i \infty\} \times \Ac_1} + \sum_{\ell \in \Hom(\Lambda_1, \Lambda_g)\setminus\{0\}}e^{\pi i \tau \mathsf{Q}(\ell)}\delta_{Z(s_\ell)} = \mathsf{e}(\mathbb F,0)  + \frac{1}{4\pi v}[\{i \infty\}\times \Ac_g]  + \sum_{d\geq 1}q^{d}[\NL_{g,d}^+]
$$
if $g=1$, which is almost holomorphic, with holomorphic part \eqref{eq:power-series-inA1xAg}. Theorem \ref{thm:new-green} follows.

\subsection{Integrability outside of $W_2$}\label{s:integrability-outside-W2}
In this section, we prove Proposition \ref{prop:integrability-noW2}. Let \(\eta\) be a test form supported on \(\mathsf B_{i\infty}(\varepsilon_0)\times\mathcal A_g\). We need to show that
\begin{equation}\label{eq:W2-int-want}
\lim_{\varepsilon\to0}
\int_{[\mathsf K_{i\infty}(\varepsilon)^c/\Gamma_\infty]\times\mathcal A_g}
\sum_{\ell\notin W_2}
e^{\pi i\tau\mathsf Q(\ell)}
\mathsf e(\mathbb F,s_{v^{1/2}\ell})\wedge\eta
\end{equation}
converges uniformly over \(\tau \in \Omega\subset \mathscr H\) compact\footnote{By the Banach--Steinhaus theorem, convergence for each $\eta$ gives uniform convergence for $\eta$ ranging over compact subsets of $\Omega^*_{c}(\mathsf{B}_{i\infty}(\eps_0) \times \mathcal A_g)$.}. Working locally on \(\mathcal A_g\), we may suppose that \(Z\) ranges in a fixed compact subset \(K\subset\mathfrak H_g\); in the argument, the constants $c,C$ are uniform over $K$ and $\Omega$, and we indicate the precise dependency by a subscript. Since $\mathfrak q_Z$ is positive definite, \(\mathfrak q_Z(\xi)\geq c_K|\xi|^2\).

We take coefficient norms in the coordinates \((z,Z)\). By Lemma \ref{lem:simplified-norm}, we have
$$
\left|
e^{\pi i\tau\mathsf Q(\ell)}
e^{-2\pi\|s_{v^{1/2}\ell}\|^2}
\right|
=
e^{-2\pi v\left(y\mathfrak q_Z(b)+y^{-1}\mathfrak q_Z(a-xb)\right)}\leq e^{-2\pi c_Kv(y|b|^2 + y^{-1}|a-xb|^2)}\,,
$$
where $|\cdot |$ is the Euclidean norm in $\mathbb R^{2g}$ and $a=(P_1, Q_1)$, $b = (P_2, Q_2)\neq 0$ since $\ell \not \in W_2$. The curvature and covariant derivative formulas \eqref{eq:curvature}, \eqref{eq:covariant-1}, \eqref{eq:covariant-2} give
\[
\|R\|\leq C_K,\qquad
\|B_{s_{v^{1/2}\ell}}\|
\leq C_Kv\left(y|b|^2+y^{-1}|a-xb|^2\right)\,.
\]
Thus, by the expression in Lemma \(\ref{lem:ex-euler}\) for the MQ Euler class, the remaining terms in the expression for
$$
e^{\pi i\tau\mathsf Q(\ell)}
\mathsf e(\mathbb F,s_{v^{1/2}\ell})
$$
are bounded by $C_K\left(1+v\left(y|b|^2+y^{-1}|a-xb|^2\right)\right)^g$,
so these factors can be absorbed into the Gaussian. Consequently,
\[
\left\|
e^{\pi i\tau\mathsf Q(\ell)}
\mathsf e(\mathbb F,s_{v^{1/2}\ell})
\right\|
\leq C_K e^{-c_K'v\left(y|b|^2+y^{-1}|a-xb|^2\right)}\,.
\]
By \eqref{eq:eta-decay}, the coefficients of $\eta$ are uniformly bounded in $z$. Therefore,
\begin{align*}
\sum_{\ell\notin W_2}\left\|\,
e^{\pi i\tau\mathsf Q(\ell)}
\mathsf e(\mathbb F,s_{v^{1/2}\ell})\wedge \eta\,\right\|&\leq C_K \sum_{b \in \mathbb Z^{2g}\smallsetminus \{0\}}e^{-c'_K vy|b|^2}\sum_{a \in \mathbb Z^{2g}} e^{-c'_K v y^{-1}|a-xb|^2}\\
&\leq C'_{K,\Omega} y^g\sum_{b\in \mathbb Z^{2g}\smallsetminus \{0\}}e^{-c'_K v y |b|^2}\\
&\leq C''_{K, \Omega} e^{-c''_{K,\Omega}y}
\end{align*}
where in the second inequality we have used that
\begin{equation}\label{eq:lattice-bound}
\sup_{h\in\mathbb R^{2g}}
\sum_{a\in\mathbb Z^{2g}}
e^{-c'_K(v/y)|a-h|^2}
=\sup_{h \in \RR^{2g}} \prod_{j=1}^{2g}\left(\sum_{n \in \ZZ}e^{-c'_K (v/y)(n-h_j)^2}\right)\leq  \left(2 + \sqrt{\frac{\pi y}{c'_Kv}}\right)^{2g} \leq C_{K,\Omega} y^{g}
\end{equation}
because $v = \mathrm{Im}(\tau)$ is bounded below. Therefore \eqref{eq:W2-int-want} converges.

\subsection{Integrability on $W_2$}\label{s:integrability-W2}
In this section, we prove Proposition \ref{prop:integrability-W2}. We need to show that the limit
\begin{equation}\label{eq:W_2-integrability-step1}
\lim_{\eps \to 0}\int_{[\mathsf{K}_{i\infty}(\eps)^c/\Gamma_\infty]\times \Ac_g}\bigg(\sum_{\ell \in W_2} \mathsf{e}(\mathbb F, s_{v^{1/2}\ell})\bigg)\wedge \eta \,,
\end{equation}
exists and converges, absolutely and uniformly in $\tau \in \mathscr{H}$. Outside of the cusp, we use Lemma \ref{lem:homogeneous-euler-W2} to write the Euler classes as
\begin{align*}
\mathsf{e}(\mathbb F, s_{v^{1/2}\ell}) =\,& \sum_{i, j \in \{0,1\}}e^{-2\pi \mathfrak q_Z((v/y)^{1/2}\ell)} E_\lambda^{(i,j)}((v/y)^{1/2}\ell) \wedge \left(\frac{dz}{y}\right)^i\wedge \left(\frac{d\overline{z}}{y}\right)^j\\
=\,& \sum_{i, j \in \{0,1\}} f^{(i,j)}((v/y)^{1/2}\ell)\wedge \left(\frac{dz}{y}\right)^i\wedge \left(\frac{d\overline{z}}{y}\right)^j\,,
\end{align*}
where $f^{(i,j)}(\ell) = e^{-2\pi \mathfrak q_Z(\ell)} E_\lambda^{(i,j)}(\ell)$ is a Schwartz function. But $f^{(i,j)}$ does not have negative exponential decay as $y \to \infty$. To remedy that, we apply Poisson summation:
$$
\sum_{\ell \in W_2} \mathsf{e}(\mathbb F, s_{v^{1/2}\ell})
= \sum_{i,j \in \{0,1\}}\bigg((y/v)^g \sum_{\xi \in W_2^\vee} \widehat{f^{(i,j)}}\Big((y/v)^{1/2}\xi\Big)\bigg)\left(\frac{dz}{y}\right)^i\wedge \left(\frac{d\overline{z}}{y}\right)^j\,,
$$
Note that the values of the Fourier transforms of the $f^{(i,j)}$ at frequency zero vanish for all $i,j$ by Proposition \ref{p:W2exact0}:
$$
\sum_{i,j} \widehat{f^{(i,j)}}(0) \left(\frac{dz}{y}\right)^i\wedge \left(\frac{d\overline{z}}{y}\right)^j = \int_{W_{2,\RR}}\mathsf{e}(\mathbb F, s_\ell)\, \mathrm{d}\ell = [\lambda^0](y\lambda)^g =0\,.
$$
If $\mathfrak q_Z^{-1}$ denotes the inverse quadratic form of $\mathfrak q_Z$ (given by the inverse of the matrix \eqref{eq:factorize-quadratic}), then
$$
\widehat{f^{(i,j)}}(\xi) = [\lambda^0]E^{(i,j)}_\lambda\Big(\frac{i}{2\pi}\partial_\xi\Big)\, e^{-2\pi \mathfrak q_Z^\vee(\xi)}\,.
$$
For $\xi \neq 0$, the corresponding term in the Poisson sum decays exponentially as $y \to \infty$. Therefore,
$$
\sum_{\ell \in W_2} \mathsf{e}(\mathbb F, s_{v^{1/2}\ell})\wedge \eta  = O(e^{-Cy/v})\,.
$$
Moreover, since the coefficients of $\mathfrak q_Z$ depend continuously on $Z$, these bounds are uniform on compact subsets of $[\mathsf{K}_{i\infty}(\eps)^c/\Gamma_\infty]\times \Ac_g$. This proves Proposition \ref{prop:integrability-W2}.

\subsection{The Green current outside of $W_2$}\label{s:gree-out-W2}

In this section, we prove Proposition \ref{prop:exact-not-W2-green}.
Let $\eta$ be a test form supported on
$\mathsf B_{i\infty}(\varepsilon_0)\times\mathcal A_g$, and let
$\Omega\subset\mathscr H$ be compact. Working locally on $\mathcal A_g$,
we may suppose that $Z$ ranges in a fixed compact subset
$K\subset\mathfrak H_g$. We take coefficient norms in the coordinates
$(z,Z)$ and assume throughout that $y\geq1$, $0\leq x\leq1$.

We first prove convergence. Unlike the Euler forms, the Green forms
are singular along $Z(s_\ell)$, so we estimate their coefficients
after integration in $Z$. By the explicit expression \eqref{eq:green-explicit} for the Green current $\mathfrak g(\mathbb F, s_\ell)$, there is an integer $N\geq0$ such
that, outside $Z(s_\ell)$,
\begin{equation}\label{eq:green-outside-pointwise}
\left\|
e^{\pi i\tau\mathsf Q(\ell)}
\mathrm d^c\mathfrak g(\mathbb F,s_{v^{1/2}\ell})
\right\|
\leq 
C_{K,\Omega}y^N(1+|\ell|)^N
\left(1+|s_\ell|^{1-2g}\right)\cdot
e^{-2\pi v\left(\|s_\ell\|^2+\frac12\mathsf Q(\ell)\right)}.
\end{equation}
Here $|s_\ell|$ denotes the Euclidean norm of the coordinate vector
of $s_\ell \in \mathbb R^{g}$.

Indeed, the term with $k=0$ in \eqref{eq:green-explicit} is a multiple of the $(g-1)$--th Chern class of the curvature, so it is $\mathrm d^c$--closed by Chern--Weil theory. For the remainder terms, one uses the exact evaluation 
\[
\int_1^\infty t^{k-1} e^{-At}\,dt
=
e^{-A}\sum_{j=0}^{k-1}
\frac{(k-1)!}{(k-1-j)!A^{j+1}},
\quad k\geq1,
\]
to obtain \eqref{eq:green-outside-pointwise}. The strongest
singularity is $|s_\ell|^{1-2g}$, and all remaining factors are
bounded by polynomials in $y$ and $|\ell|$, uniformly for
$Z\in K$ and $\tau\in\Omega$.

We can assume that $K$ is a Euclidean ball
$$
B_R = \Big \{Z \in \mathrm{Mat}_{g\times g}(\mathbb C) : Z= Z^t\,,  \sum_{i, j} |Z_{ij}|^2 \leq R\Big\}\,.
$$
Recall that  $s_\ell(z,Z) = Z(zP_2-P_1) + zQ_2-Q_1$. If $zP_2-P_1=0$ then, since $z$ has purely imaginary part, $P_1=P_2=0$, so $Q_2\neq 0$ because $\ell \not \in W_2$, so $|s_\ell(z,Z)| = |zQ_2-Q_1|\geq y|Q_2|\geq 1$ for all $\ell$. If $zP_2-P_1\neq0$, choose a unitary $U$ with
$$
Ue_1=\frac{zP_2-P_1}{|zP_2 - P_1|}
$$
Let $W=U^T ZU$. Since $U^T$ is unitary,
\[
|s_\ell(z,Z)|=|zP_2 - P_1|\bigg|We_1+\frac{U^T(zQ_2-Q_1)}{|zP_2-P_1|}\bigg|.
\]
And note that $\frac{U^T(zQ_2-Q_1)}{|zP_2-P_1|}$ is independent of $Z$ and $|zP_2 - P_1|\geq 1$. We change variables by $Z \mapsto U^TZU$ (which preserves the volume form and $B_R$) and integrate over the entries of $W$ outside its first column. This gives
\begin{align*}
\int_{B_R}  |s_\ell(z,Z)|^{1-2g} dV_Z &= |zP_2-P_1|^{1-2g} \int_{B_R} \bigg|We_1+\frac{U^T(zQ_2-Q_1)}{|zP_2-P_1|}\bigg|^{1-2g} dV_W\\
&\leq C_R \int_{\{w \in \mathbb C^g : |w|\leq R^2\}} \bigg|w + \frac{U^T(zQ_2-Q_1)}{|zP_2-P_1|}\bigg|^{1-2g} dw\\
&\leq C_R\bigg( \mathrm{vol} (\{w \in \mathbb C^g : |w|\leq R^2\}) +  \int_{|w|\leq R^2} |w|^{1-2g}dw\bigg)\,.
\end{align*}
In both cases, we conclude that
\begin{equation}\label{eq:green-outside-singularity}
\int_K\left(1+|s_\ell(z,Z)|^{1-2g}\right)\,dV_Z
\leq C_K
\end{equation}
uniformly in $z$ and $\ell\notin W_2$.

Write $a=(P_1,Q_1)$ and $b=(P_2,Q_2)\neq0$. By
Lemma \ref{lem:simplified-norm},
\[
\|s_\ell\|^2+\frac12\mathsf Q(\ell)
=
y\mathfrak q_Z(b)+y^{-1}\mathfrak q_Z(a-xb)
\geq c_K\left(y+y^{-1}|\ell|^2\right).
\]
Here we use $|b|\geq1$, $y\geq1$, and
$|a|^2\leq2|a-xb|^2+2|b|^2$.
Combining this with
\eqref{eq:green-outside-pointwise} and
\eqref{eq:green-outside-singularity}, we obtain
\begin{align*}
\sum_{\ell\notin W_2}\int_K
\left\|
e^{\pi i\tau\mathsf Q(\ell)}
\mathrm d^c\mathfrak g(\mathbb F,s_{v^{1/2}\ell})
\right\|\,dV_Z
&\leq
C_{K,\Omega}y^N e^{-2\pi c_Kvy}
\sum_{\ell\in\mathbb Z^{4g}}
(1+|\ell|)^N e^{-2\pi c_K(v/y)|\ell|^2}\\
&\leq
C'_{K,\Omega}y^{2g+3N/2}e^{-2\pi c_Kvy}\\
&\leq
C''_{K,\Omega}e^{-c'_{K,\Omega}y}.
\end{align*}
The second estimate follows by a similar lattice sum as in \eqref{eq:lattice-bound} and the last inequality uses that $v$ is bounded below in $\Omega$. By \eqref{eq:eta-decay}, the coefficients of $\eta$ are bounded.
Integrating the preceding estimate over 
$y\geq1$ proves absolute convergence, uniformly for
$\tau\in\Omega$. We denote the resulting current by
$\overline{\Psi}_{W_2^c}(\tau)$.

% Put $u=zP_2-P_1$, so that
% \[
% s_\ell(z,Z)=Zu+zQ_2-Q_1.
% \]
% If $u\neq0$, integrality of $P_1,P_2$ and $y\geq1$ imply
% $|u|\geq1$. Choose $j$ such that
% $|u_j|\geq |u|/\sqrt g$. Keeping the other entries of $Z$ fixed,
% the affine map from the $g$ independent entries in the $j$-th
% row and column of $Z$ to $s_\ell$ has complex Jacobian
% determinant $u_j^g$, and its linear inverse has norm
% $O(|u|^{-1})$. Integration in these $g$ complex variables gives
% \[
% \int_K |s_\ell|^{1-2g}\,dV_Z
% \leq C_K|u|^{1-2g}\leq C_K.
% \]
% The estimate is uniform in the translation term, since
% $|w|^{1-2g}$ is locally integrable on $\mathbb C^g$.
% If $u=0$, then $P_1=P_2=0$, and
% \[
% |s_\ell|=|zQ_2-Q_1|\geq y|Q_2|\geq1,
% \]
% because $\ell\notin W_2$. This proves
% \eqref{eq:green-outside-singularity}.

It remains to prove the Green equation \eqref{eq:Green-outside-W2} across the cusp. Let
$$
\omega(\tau) = \sum_{\ell \in \Hom(\Lambda_1, \Lambda_g)\smallsetminus W_2}e^{\pi i \tau \mathsf{Q}(\ell)}\mathsf{e}(\mathbb F, s_{v^{1/2}\ell}) - \sum_{\ell \in \Hom(\Lambda_1, \Lambda_g)\smallsetminus \{0\}} e^{\pi i \tau \mathsf{Q}(\ell)}\delta_{Z(s_\ell)}\,.
$$
Choose smooth
functions $\chi_T(y)$ which are equal to $1$ for $y\leq T$,
equal to $0$ for $y\geq T+1$, and satisfy
$|\d\chi_T|\leq C$ independently of $T$.
Then, for every test form $\eta$, we have
\begin{align*}
\int_{\mathsf{B}_{i\infty}(\eps_0)\times \mathcal A_g} \d \overline{\Psi}_{W_2^c}(\tau)\wedge \eta &= \int_{\mathsf{B}_{i\infty}(\eps_0)\times \mathcal A_g}  \overline{\Psi}_{W_2^c}(\tau)\wedge \d\eta  \\
&=\lim_{T\to \infty} \int_{\mathsf{B}_{i\infty}(\eps_0)\times \mathcal A_g}  \overline{\Psi}_{W_2^c}(\tau)\wedge \chi_T\d\eta - \overline{\Psi}_{W_2^c}(\tau)\wedge \d\chi_T \wedge \eta\\
&=\lim_{T\to \infty} \int_{\mathsf{B}_{i\infty}(\eps_0)\times \mathcal A_g}\chi_T \omega(\tau)\wedge \eta - \lim_{T\to \infty} \int_{[\mathsf{B}_{i\infty}(\eps_0)]\times \mathcal A_g} \overline{\Psi}_{W_2^c}(\tau)\wedge \d\chi_T\wedge \eta\\
&= \int_{\mathsf{B}_{i\infty}(\eps_0)\times \mathcal A_g} \omega(\tau)\wedge \eta - \lim_{T\to \infty} \int_{[\mathsf{B}_{i\infty}(\eps_0)]\times \mathcal A_g} \overline{\Psi}_{W_2^c}(\tau)\wedge \d\chi_T\wedge \eta
\end{align*}
where in the third equality we have used that by Lemmas \ref{lem:reg-section} and \ref{l:g-curr-dZ}, $\d \overline{\Psi}_{W_2^c}(\tau) = \omega(\tau)$ holds on the support of $\chi_T$, and all commutations of limits are allowed because we already know that $\overline{\Psi}_{W_2^c}(\tau)$ is a current, and by Lemma \ref{lem:convergence-deltas} and Proposition \ref{prop:integrability-noW2}, $\omega(\tau)$ is also a current on $\mathsf{B}_{i\infty}(\eps_0)\times \mathcal A_g$. By \eqref{eq:green-outside-pointwise} and the assumption on $\d\chi_T$, the last term in the above derivation tends to $0$:
\[
\left|
\int
\sum_{\ell\notin W_2}
e^{\pi i\tau\mathsf Q(\ell)}
\mathrm d^c\mathfrak g(\mathbb F,s_{v^{1/2}\ell})
\wedge\mathrm d\chi_T\wedge\eta
\right|
\leq
C_{K,\Omega}\int_T^{T+1}e^{-c'_{K,\Omega}y}\,dy
\longrightarrow0.
\]
Therefore, $\d\overline{\Psi}_{W_2^c}(\tau) = \omega(\tau)$, which finishes the proof of Proposition \ref{prop:exact-not-W2-green}.

\subsection{The singularities of the Green current on $W_2$}\label{s:green-in-W2}
In this section, we prove Proposition \ref{prop:exact-W2-green}. We need to study the sum
\begin{equation}\label{eq:W_2-green-step1}
\int_{\mathsf{B}_{i\infty}(\eps_0)\times \Ac_g} \sum_{\ell \in W_2\smallsetminus \{0\}} \mathsf{e}(\mathbb F, s_{v^{1/2}\ell})\wedge \eta \,,
\end{equation}
where $\eta$ satisfies \ref{notation} and $\operatorname{supp}(\eta)\subset \mathsf{B}_{i\infty}(\eps_0) \cong [\mathsf{K}_{i\infty}(\eps_0)/\Gamma_\infty]\cup \{i \infty\}$. By Theorem \ref{thm:integrability}, we can evaluate the above integral by taking a neighborhood $\mathsf{B}_{i\infty}(\eps)$ around the cusp $i \infty \in \Abar_1$ and instead calculating
\begin{equation}
    \lim_{\eps \to 0} \Phi(\eps, \eta)\,, \quad \Phi(\eps, \eta):= \int_{\Ac_g \times [\mathsf{K}_{i\infty}(\eps)^c/\Gamma_\infty]}\sum_{\ell \in W_2 \smallsetminus \{0\}}\mathsf{e}(\mathbb F, s_{v^{1/2}\ell})\wedge \eta\,.
\end{equation}
\begin{proof}[Proof of Proposition \ref{prop:exact-W2-green}] We show integrability first.
By the definition of the Green form and Lemma \ref{lem:homogeneous-euler-W2}, we have
\begin{align*}
\mathrm{d}^c\mathfrak g(\mathbb F, s_{v^{1/2}\ell}) =\,& [\lambda^1] \int_{1}^\infty\mathrm{d}^c\mathsf{e}_{\lambda}(\mathbb F, (tv)^{1/2}s_{\ell})\frac{dt}{t}\\
=\,& \sum_{i, j \in \{0,1\}}[\lambda^1]\left(\int_{1}^\infty e^{-2\pi \mathfrak q_Z((tv/y)^{1/2}\ell)}F_\lambda^{(i,j)}((tv/y)^{1/2}\ell)\frac{dt}{t}\right)\left(\frac{dz}{y}\right)^i\wedge \left(\frac{d\overline{z}}{y}\right)^j \\
=\,& \sum_{i, j \in \{0,1\}}[\lambda^1]\left(\int_{v/y}^\infty e^{-2\pi \mathfrak q_Z(t^{1/2}\ell)}F_\lambda^{(i,j)}(t^{1/2}\ell)\frac{dt}{t}\right)\left(\frac{dz}{y}\right)^i\wedge \left(\frac{d\overline{z}}{y}\right)^j\,,%\\
%=\,&\sum_{i, j \in \{0,1\}}[\lambda^1]\left(\int_{v/y}^\infty e^{-2\pi \mathfrak q_Z(t^{1/2}\ell)}E_\lambda^{(i,j)}(t^{1/2}\ell)\frac{dt}{t}\right)\left(\frac{dz}{y}\right)^i\wedge \left(\frac{d\overline{z}}{y}\right)^j
\end{align*}
where $F^{(i,j)}_{\lambda}$ are polynomials valued in $A^{*-i, *-j}(\mathfrak H_g)[\lambda]$. The exchange of $\mathrm{d}^c$ and the integral is justified because $s_\ell$ does not vanish, so both $\mathsf{e}_\lambda(\mathbb F, (tv)^{1/2}s_\ell)$ and its derivatives are uniformly bounded by $e^{-Ctv}$. We want to apply Lemma \ref{lem:analysis} below, with $f(\xi)= f^{(i,j)}(\xi) = e^{-2\pi\mathfrak q_Z(\xi)}F_{\lambda}^{(i,j)}(\xi)$ for $i,j \in \{0,1\}$. The constant terms vanish.
$$
\sum_{i,j \in \{0,1\}}F_{\lambda}^{(i,j)}(0) \left(\frac{dz}{y}\right)^i\wedge \left(\frac{d\overline{z}}{y}\right)^j = \mathrm{d}^c \mathsf{e}_\lambda (\mathbb F, 0) = 0\,.
$$
The last equality holds because the Chern--Weil form $\mathsf{e}_\lambda(\mathbb F,0)$ is $\mathrm{d}^c$-closed. By Proposition \ref{p:W2exact0}, we have
\begin{equation}\label{eq:green-fourier-zero}
\sum_{i,j \in \{0,1\}}\widehat{f^{(i,j)}}(0) \left(\frac{dz}{y}\right)^i\wedge \left(\frac{d\overline{z}}{y}\right)^j= [\lambda^1]\mathrm{d}^c \int_{W_{2,\RR}}\mathsf{e}_\lambda (\mathbb F, s_\ell)\, \mathrm{d}\ell = \begin{cases}
        -\frac{dz + d\overline{z}}{8\pi y}&\text{ if } g=1\,,\\
        0&\text{ otherwise.}
    \end{cases}
\end{equation}
Therefore, by Lemma \ref{lem:analysis}, we have
$$
\overline{\Psi}_{W_2}(\tau) = \sum_{\ell \in W_2 \smallsetminus \{0\}}\mathrm{d}^c \mathfrak g(\mathbb F, s_{v^{1/2}\ell}) = O(1) + O(y^{-1})(dz + d\overline{z}) + O(y^{-2})dz\wedge d\overline{z} + O(e^{-Cy/v})\,.
$$
The statement is independent of $g$, since the $y$ in the denominator of $\frac{dz+d\overline{z}}{8\pi y}$ cancels with $y^g$ from Lemma \ref{lem:analysis} when $g=1$. This shows integrability.

By Lemma \ref{lem:reg-section}, $s_{\ell}$ has no zeros for all $\ell \in W_2 \smallsetminus \{0\}$, so we can apply the second part of Lemma \ref{l:g-curr-dZ} outside of the cusp. By Stokes' theorem, we have
\begin{align}
\begin{split}\label{eq:Stokes}
     \Phi(\eps, \eta) =\,& \int_{[\mathsf{K}_{i\infty}(\eps)^c/\Gamma_{\infty}] \times \Ac_g}\sum_{\ell \in W_2 \smallsetminus \{0\}}\ddc \mathfrak g (\mathbb F, s_{v^{1/2}\ell})\wedge \eta \\
    =\, & \int_{[\mathsf{K}_{i\infty}(\eps)^c/\Gamma_\infty] \times \Ac_g} \sum_{\ell \in W_2 \smallsetminus \{0\}}\mathrm{d}^c\mathfrak g (\mathbb F, s_{v^{1/2}\ell})\wedge d\eta + \int_{\mathsf{S}_{i\infty}(\eps) \times \Ac_g}\sum_{\ell \in W_2 \smallsetminus \{0\}}\mathrm{d}^{c}\mathfrak g (\mathbb F, s_{v^{1/2}\ell})\wedge \eta
\end{split}
\end{align}
where $\mathsf{S}_{i\infty}(\eps) = \partial \mathsf{K}_{i\infty}(\eps) \cong \partial \mathsf{B}_{i\infty}(\eps)$ is the circle $|p| = \eps$, oriented negatively.

Since $\eta$ is smooth, the first term in \eqref{eq:Stokes} converges to
$$
\lim_{\eps\to 0} \int_{[\mathsf{K}_{i\infty}(\eps)^c/\Gamma_{\infty}]\times \Ac_g} \sum_{\ell \in W_2 \smallsetminus \{0\}}\mathrm{d}^c\mathfrak g (\mathbb F, s_{v^{1/2}\ell})\wedge \mathrm{d}\eta = \int_{[\mathsf{K}_{i\infty}(\eps_0)/\Gamma_\infty]\times \Ac_g}\sum_{\ell \in W_2\smallsetminus \{0\}} \mathrm{d}^c\mathfrak g(\mathbb F, s_{v^{1/2}\ell})\wedge \mathrm{d} \eta
$$
which is precisely the definition of $\d\overline{\Psi}_{W_2}(\tau) \wedge \eta$. Now we analyze the second term.

We decompose the form $\eta$ into its $(i,j)$ components as in \ref{notation} and apply Lemma \ref{lem:analysis}. After restricting to the sphere $\mathsf{S}_{i\infty}(\eps)$ and using Lemma \ref{lem:analysis} as before (with $f =f^{(i,j)}= e^{-2\pi\mathfrak q_Z(\xi)}F_{\lambda}^{(i,j)}(\xi)$), and \eqref{eq:green-fourier-zero} we see that the second term in \eqref{eq:Stokes} is given by
\begin{align*}
\sum_{\ell \in W_2 \smallsetminus \{0\}}\mathrm{d}^{c}\mathfrak g (\mathbb F, s_{v^{1/2}\ell})\wedge \eta\,\Big|_{\mathsf{S}_{i\infty}(\eps) \times \Ac_g} =\,& K(f^{(0,0)})(\eta^{(1,0)}dz + \eta^{(0,1)}d\overline{z}) + \eta^{(0,0)}\frac{K(f^{(1,0)})dz + K(f^{(0,1)})d\overline{z}}{y}\\
& - \begin{cases}
    \eta^{(0,0)}\frac{dz + d\overline{z}}{8v\pi}&\text{ if }g=1\\
    0&\text{otherwise}
\end{cases} + O(e^{-Cy/v})\,,
\end{align*}
where $C, K(f^{(i,j)})$ are uniform in $Z$. By the estimates \eqref{eq:eta-decay} on $\eta^{(i,j)}$, we find that
\begin{equation}
    \sum_{\ell \in W_2 \smallsetminus \{0\}}\mathrm{d}^{c}\mathfrak g (\mathbb F, s_{v^{1/2}\ell})\wedge \eta\,\Big|_{\mathsf{S}_{i\infty}(\eps) \times \Ac_g} =  \begin{cases}
    - \eta^{(0,0)}\frac{dz + d\overline{z}}{8v\pi}&\text{ if }g=1\\
    0&\text{otherwise}
\end{cases} + O(1/y)
\end{equation}
For $z = x+iy \in \mathsf{S}_{i\infty}(\eps)$ we have $y = -\frac{\log \eps}{2\pi} \to \infty$, so
\begin{align*}
    \lim_{\eps \to 0} \int_{\mathsf{S}_{i\infty}(\eps) \times \Ac_g}\sum_{\ell \in W_2 \smallsetminus \{0\}}\mathrm{d}^{c}\mathfrak g (\mathbb F, s_{v^{1/2}\ell})\wedge \eta=\,& \lim_{\eps \to 0}\int_{\mathsf{S}_{i\infty}(\eps)\times \Ac_g} \eta^{(0,0)}\wedge  \begin{cases}
    - \frac{dz + d\overline{z}}{8v\pi}&\text{ if }g=1\\
    0&\text{otherwise}
    \end{cases}
\end{align*}
This gives $\int \frac{1}{4\pi v}\delta_{\{i \infty\} \times \Ac_g}\wedge \eta$, since $\mathsf{S}_{i\infty}(\eps)$ is negatively oriented.
\end{proof}

\begin{lemma}\label{lem:analysis}
    Let $V$ be a real vector space of dimension $2g$ with a lattice $L$ of covolume $1$, and let $L^\vee$ be the dual lattice. Let $f(x) = p(x)e^{-2\pi\mathfrak q(x)}$, where $p$ is a polynomial and $\mathfrak q$ is a positive definite quadratic form on $V_{\RR}$, with $f(0) = 0$, and let $\widehat f(\xi) = \int_V f(x)e^{-2\pi i\langle x,\xi\rangle}dx$ be its Fourier transform. We write
    $$
    S(\sigma) = \sum_{\ell\in L}f(\sigma^{1/2}\ell)\,,\quad R(\sigma) = \sigma^{-g}\sum_{\ell'\in L^\vee\smallsetminus\{0\}}\widehat f(\sigma^{-1/2}\ell')\,,\quad K(f) = -\frac{\widehat f(0)}{g} + \int_0^1R(\sigma)\frac{d\sigma}{\sigma} + \int_1^\infty S(\sigma)\frac{d\sigma}{\sigma}\,,
    $$
    Then, as $a$ approaches $0$, we have\footnote{If $f$ is a Schwartz function, we get the same equation, but with an error term that is rapidly decaying.}
    $$
    \sum_{\ell\in L}\int_a^\infty f(t^{1/2}\ell)\,\frac{dt}{t} = \frac{\widehat f(0)}{g}\,a^{-g} + K(f) + O(e^{-C/a})\,,
    $$
    The constant $C$ and the bound on $K(f)$ are uniform when $f$ depends continuously on a parameter $Z$ ranging in a compact set.
\end{lemma}

\begin{proof}
    For $\ell\neq0$ and $t\geq a$ we have
    $$
    |f(t^{1/2}\ell)|\leq c(1+t|\ell|^2)^Ne^{-2\pi C_0t|\ell|^2}\,,
    $$
    where $C_0$ is the smallest eigenvalue of $\mathfrak q$. This is summable over $\ell$ and integrable in $\frac{dt}{t}$, so by dominated convergence, we have
    $$
    \sum_{\ell\in L}\int_a^\infty f(t^{1/2}\ell)\,\frac{dt}{t} = \int_a^{\infty}S(\sigma)\frac{d\sigma}{\sigma}\,.
    $$ 

    The function $x\mapsto f(\sigma^{1/2}x)$ has Fourier transform $\sigma^{-g}\widehat f(\sigma^{-1/2}\xi)$. By Poisson summation,
    $$
    S(\sigma) = \sigma^{-g}\widehat f(0) + R(\sigma)\,.
    $$
    Splitting the integral at $\sigma = 1$ yields
    $$
    \int_a^\infty S(\sigma)\frac{d\sigma}{\sigma} = \int_1^\infty S(\sigma)\frac{d\sigma}{\sigma} + \widehat f(0)\int_a^1\sigma^{-g-1}d\sigma + \int_0^1R(\sigma)\frac{d\sigma}{\sigma} - \int_0^aR(\sigma)\frac{d\sigma}{\sigma}\,.
    $$
    This is the claimed identity with the error term controlled by $ -\int_0^aR(\sigma)\frac{d\sigma}{\sigma}$. The Fourier transform of $f$ is also a negative Gaussian times a polynomial, so
    $$
    |R(\sigma)|\leq c\sigma^{-g}\sum_{\ell'\neq0}\Big(1+\frac{|\ell'|}{\sigma^{1/2}}\Big)^Ne^{-C'|\ell'|^2/\sigma}= O(e^{-C'/\sigma})\,.
    $$
    Every constant above depends continuously on the coefficients of $f$, so the last claim follows.
\end{proof}

\section{Hecke operators}\label{sec:hecke}

Recall from \cite{CMOP24} that there is a natural \emph{tautological projection} operator, taking any pure weight class to a polynomial in the $\lambda_i$ classes:
$$
\taut : H_{\mathsf{pure}}^*(\Ac_g) \to RH^*(\Ac_g)\,.
$$
Since $\taut \circ \taut = \taut$, it provides a canonical splitting
$$
(\taut, \operatorname{id} - \taut) : H^*_{\mathsf{pure}}(\Ac_g) \to RH^*(\Ac_g) \oplus H_{\mathsf{pure}}^*(\Ac_g)^{\mathsf{cusp}}\,,
$$
where $H_{\mathsf{pure}}^*(\Ac_g)^{\mathsf{cusp}}:= \ker (\taut)$\footnote{The Hecke algebra of $\mathcal A_g$ decomposes the cohomology $H^*(\mathcal A_g)$ into common eigenspaces of the Hecke operators. In different contexts, $H^{*}_{cusp}(\mathcal A_g)$ denotes the part whose eigenvalue is not the volume. These two notions coincide, since volume eigenvalues are precisely the tautological ring by \cite[Proposition 2.3]{EsnaultHarris2019}.}.

As shown in \cite{IribarLopez}, $\taut (\mathbf{NL}_{g}) = \QQ\lambda_{g-1}$, so the above splitting induces a splitting
$$
\mathbf{NL}_{g} = \QQ\lambda_{g-1} \oplus \mathbf{NL}_g^{\mathsf{cusp}}\,,
$$
where $\mathbf{NL}_g^{\mathsf{cusp}}$ is spanned by the cycles
$$
[\NL_{g,d}] - \taut ([\NL_{g,d}]) = [\NL_{g,d}] - \frac{g}{6|B_{2g}|}\sigma_{2g-1}(d)\lambda_{g-1}\,.
$$

\begin{lemma}\label{lem:Psi-cusp}
    Consider the linear map $\psi : \operatorname{Mod}_{2g}^\vee \to \mathbf{NL}_g$ from \eqref{eq:Psi-map}. Let $(\operatorname{Mod}_{2g})^\vee \to (\operatorname{SMod}_{2g})^\vee$ be the restriction map. Then, there is a natural surjective map $\psi^{\mathsf{cusp}}$ that makes the following diagram commute:
    $$
    \begin{tikzcd}
	{(\operatorname{Mod}_{2g})^\vee} & {\mathbf{NL}_g} \\
	{(\operatorname{SMod}_{2g})^\vee} & {\mathbf{NL}_g^{\mathsf{cusp}}}
	\arrow["\psi", from=1-1, to=1-2]
	\arrow[from=1-1, to=2-1]
	\arrow["{\operatorname{id} - \taut}", from=1-2, to=2-2]
	\arrow["{\psi^{\mathsf{cusp}}}"', from=2-1, to=2-2]
\end{tikzcd}
    $$
\end{lemma}
\begin{proof}
    The kernel of $(\operatorname{Mod}_{2g})^\vee \to (\operatorname{SMod}_{2g})^\vee$ is generated by $a_0(-)$, the functional that extracts the constant Fourier coefficient, and we have
    $$
    \psi(a_0(-)) - \taut(\psi(a_0(-))) = \frac{(-1)^g}{24} \big(\lambda_{g-1} - \taut(\lambda_{g-1})\big) =0
    $$
\end{proof}

\begin{remark}\label{rem:NL-taut}
    From Lemma \ref{lem:Psi-cusp}, it is clear that the subspace of $H^{2g-2}(\mathcal A_g)$ spanned by tautological classes and $\mathrm{NL}$ cycles can be written as $\mathrm{RH}^{2g-2}(\mathcal A_g) \oplus \mathbf{NL}_g^{\mathsf{cusp}}$, with
    $$
    \mathrm{dim}(\mathbf{NL}_g^{\mathsf{cusp}}) \leq \dim (\mathrm{SMod}_{2g})\,,
    $$
    and since $\lambda_{g-1}\neq 0$, equality holds if and only if $\psi$ is an isomorphism.
\end{remark}

We consider geometric Hecke operators acting on the cohomology of $\Ac_g$. The construction of these operators is classical using double cosets, but we briefly review it in the simplest case.
\begin{definition}
    Let $\mathcal B_{g,n}$ parametrize triples $((A, \theta_A), (B, \theta_B), f)$, where $(A,\theta_A)$ and $(B, \theta_B)$ are principally polarized abelian varieties of dimension $g$, $f : A \to B$ is an isogeny and $f^* \theta_B = n \theta_A$.
    
    The forgetful maps $p_A, p_B : \mathcal B_{g,n} \to \mathcal A_g$ are finite \'etale maps. The \emph{geometric} Hecke operator is
    $$
    H_n := (p_B)_*(p_A)^* : H^*(\Ac_g) \to H^*(\Ac_g)\,.
    $$
    Its volume is given by $\operatorname{vol}(H_n) = H_n(1) = \deg(p_B) \in \ZZ$. Let
    $$
    H_n^{\mathsf{cusp}} = \frac{1}{\operatorname{vol}(H_n)} H_n - \operatorname{id}\,.
    $$
\end{definition}
It is clear that for any tautological class $\alpha$, we have $H_n(\alpha) = \operatorname{vol}(H_n)\alpha$. Therefore, the Hecke operators $H_n^{\mathsf{cusp}}$ act naturally on $H^*_{\mathsf{pure}}(\Ac_g)^{\mathsf{cusp}}$. As noted in \cite{IribarLopez}, the Hecke operators also preserve $\mathbf{NL}_g$. Therefore, $H_n^{\mathsf{cusp}}$ acts on $\mathbf{NL}_{g}^{\mathsf{cusp}}$.

Classically, Hecke operators act on the space of modular forms by similar correspondences
$$
T_n : \operatorname{Mod}_{2g} \to \operatorname{Mod}_{2g}\,.
$$
The action on $(\operatorname{Mod}_{2g})^\vee$ is given by
\begin{equation}\label{eq:Hecke-action}
T_n^\vee (a_d(-)) = \sum_{m \mid \gcd (d,n)}m^{2g-1}a_{nd/m^2}(-)
\end{equation}
Then, $T_n^\vee (a_0(-)) = \sigma_{2g-1}(n)a_0(-)$, where $\sigma_{2g-1}(n)$ is the sum of the $2g-1$--th powers of the divisors of $n$. The cuspidal Hecke operators are defined to be
$$
T_p^{\mathsf{cusp}} = \frac{1}{(1+p^{g-1})(1+p^{g})}\big(T_p^\vee -\sigma_{2g-1}(p)\operatorname{id}\big)\curvearrowright (\operatorname{Mod}_{2g})^\vee
$$
It is clear that the action descends to $(\operatorname{SMod}_{2g})^\vee$.
\begin{proposition}\label{prop:compatible-Hecke}
    Let $\mathcal H(\SL_2(\ZZ))^{\mathsf{cusp}}$ be the (non-unital) commutative $\QQ$-algebra generated by symbols $H_p^{\mathsf{cusp}}$, for $p$ a prime number. It acts naturally on $\mathbf{NL}_g^{\mathsf{cusp}}$, and on $(\operatorname{SMod}_{2g})^\vee$ by $T_p^{\mathsf{cusp}}$. Then, the map $\psi^{\mathsf{cusp}}$ from Lemma \ref{lem:Psi-cusp} is a homomorphism of $\mathcal H(\SL_2(\ZZ))^{\mathsf{cusp}}$-modules.
\end{proposition}
\begin{proof}
    Throughout the proof, terms of the form $a_{d/p}(-)$ or $[\NL_{g,d/p}]$ are understood to be $0$ if $p$ does not divide $d$. By \eqref{eq:Hecke-action}, we have
    $$
    T_p^{\mathsf{cusp}} (a_d(-)) = \frac{a_{pd}(-)+p^{2g-1}a_{d/p}(-)-(1+p^{2g-1})a_d(-)}{(1+p^{g-1})(1+p^g)}
    $$
    so we need to show that
    \begin{equation}\label{eq:Hecke-desire}
    H_p [\NL_{g,d}] = \frac{\operatorname{vol}(H_p)}
{(1+p^{g-1})(1+p^g)}\left([\NL_{g,pd}]+p^{2g-1}[\NL_{g,d/p}]+(p^{g-1}+p^g)[\NL_{g,d}]\right)
    \end{equation}
    The left side of \eqref{eq:Hecke-desire} is calculated by forgetting everything except $(B, \theta_B)$ on the stack that parametrizes a sequence of homomorphisms
    $$
    E \stackrel{u} \to A \stackrel{f}\to B
    $$
    where $(A, \theta_A)$ and $(B,\theta_B)$ are principally polarized abelian varieties, $E$ is an elliptic curve and we have $u^*\theta_A = d\theta_E$, $f^* \theta_B = p \theta_A$. The maps $f: A \to B$ such that $f^*\theta_B = p\theta_A$ correspond bijectively to maximal isotropic subgroups $K \subset B[p] \cong \mathbb F_p^{\oplus 2g}$. Moreover, given $h:E \to B$, it factors through $f$ if and only if $h(E[p]) \subset K$. For $r =0,1,2$, let
    $$
    Z_r \subset \NL_{g,pd}
    $$
    be the locus where the rank of $h(E[p])\subset B[p]$ as an $\mathbb F_p$-vector space is $r$. Since there are
    $$
    \prod_{i=1}^{g-r}(1+p^i)
    $$
    maximal isotropic subspaces containing $E[p]$, we obtain
    \begin{equation}\label{eq:middle-stepp}
    H_p [\NL_{g,d}] = \prod_{i=1}^{g}(1+p^i)[Z_0] + \prod_{i=1}^{g-1}(1+p^i)[Z_1] + \prod_{i=1}^{g-2}(1+p^i) [Z_2]\,.
    \end{equation}
    
    Note that if $h(E[p]) =0$ then the map $E \to B$ factors uniquely as $E \stackrel{[p]}\to E \to B$, so $[Z_0] = [\NL_{g,d/p}]$. Consider the auxiliary stack $Z$ that parametrizes maps $h: E \to B$ of degree $pd$ and a subgroup $C \subset \ker (h)\cap E[p]$ of order $p$. Forgetting $C$ shows that $[Z] = (p+1)[Z_0] + [Z_1]$ as cycles on $\NL_{g,pd}$. On the other hand, quotienting $C$ gives a degree $d$ map $E/C \to B$, so $[Z] = (p+1)[\NL_{g,d}]$. This shows that $[Z_1] = (p+1)\big([\NL_{g,d}] - [\NL_{g,d/p}]\big)$. Since $[Z_0] + [Z_1] + [Z_2] = [\NL_{g,dp}]$, we also obtain $[Z_2] = [\NL_{g,pd}]  -(p+1)[\NL_{g,d}] + p[\NL_{g,d/p}]$.

    If we substitute these formulas into \eqref{eq:middle-stepp} and use that $\operatorname{vol}(H_p) = \prod_{i=1}^{g} (1+p^i)$, we obtain \eqref{eq:Hecke-desire}.
\end{proof}
\begin{corollary}\label{cor:kernel-is-Hecke-equiv}
    The kernel of $\psi^{\mathsf{cusp}}$ is a $\mathcal H(\SL_2(\ZZ))$-submodule of $(\operatorname{SMod})^\vee$.
\end{corollary}
\begin{proof}[Proof of Theorem \ref{thm:non-vanish}]
    By \cite{vdG99}, the class $\lambda_{g-1}$ is nonzero. Therefore,  $\psi$ also restricts to an isomorphism $\QQ a_0(-) \to \QQ\lambda_{g-1}$ between the kernels of the vertical maps in Lemma \ref{lem:Psi-cusp}, and $\psi$ is injective if and only if $\psi^{\mathsf{cusp}}$ is injective. By Corollary \ref{cor:kernel-is-Hecke-equiv}, the kernel of $\psi^{\mathsf{cusp}}$ is a Hecke-invariant subspace. If Conjecture \ref{conj:modular-forms} holds, then such a space is either empty or the total space. This implies the result.
\end{proof}
\begin{proof}[Proof of Corollary \ref{cor:psi-is-iso}]
    If $g = 2,3,4,5,7$ then $\operatorname{SMod}_{2g} =0$ so the result is clear. By \cite{COP24, CLS26}, the class $[\mathcal A_1 \times \mathcal A_{5}]$ is not tautological. By \cite{IribarLopez}, for $g =12$ or $g \geq 16$ even, the class $[\mathcal A_1 \times \mathcal A_{g-1}]$ is not tautological. Therefore, part (b) of Theorem \ref{thm:non-vanish} cannot hold.
\end{proof}

\section{Geometric applications}\label{sec:geometry}

\begin{proof}[Proof of Theorem \ref{thm:lg-1}]
    By formula \eqref{eq:low-g}, we can assume that $g \geq 6$. Applying M\"{o}bius inversion to the relationship \eqref{eq:NL-from-prim} between $[\NL_{g,d}]$ and $[\NL_{g,d}^{\mathsf{prim}}]$ gives
    \begin{equation}\label{eq:Mobius}
    [\NL_{g,d}^{\mathsf{prim}}] - \taut([\NL_{g,d}^{\mathsf{prim}}]) = \sum_{ek \mid d}\mu(e)\mu(k)k\big([\NL_{g,d/(ek)}] - \taut([\NL_{g,d/(ek)}])\big)\,.
    \end{equation}
    Since the power series $\sum_{d \geq 1} \big([\NL_{g,d}] - \taut([\NL_{g,d}])\big) q^d$ is the Fourier expansion of a cusp form, we have
    $$
    [\NL_{g,d}] - \taut([\NL_{g,d}]) = O(d^{g}).
    $$
    And by \eqref{eq:Mobius}, we also have
    $$
    [\NL_{g,d}^{\mathsf{prim}}] - \taut([\NL_{g,d}^{\mathsf{prim}}]) = O\big(d^g\sum_{e\mid d}\frac{\sigma_1(e)}{e^g}\big) = O(d^g)\,,\quad (g\geq 3)\,.
    $$
    The tautological projection of $[\NL_{g,d}^{\mathsf{prim}}]$ was calculated in \cite{IribarLopez}:
    $$
    \taut([\NL_{g,d}^{\mathsf{prim}}]) = \frac{gd^{2g-1}\prod_{p \mid d\text{ prime}}(1-p^{2-2g})}{6|B_{2g}|}\lambda_{g-1}\,.
    $$
    Since $\sigma_{2g-1}(d)\geq d^{2g-1}$ and $\prod_{p \mid d\text{ prime}} (1-p^{2-2g}) \geq \prod_{p \text{ prime}} (1-p^{2-2g}) = \frac{1}{\zeta(2g-2)}>0$, we obtain
    \begin{equation}\label{eq:limits}
    \lim_{d \to \infty}\frac{6|B_{2g}|[\NL_{g,d}]}{g\sigma_{2g-1}(d)} = \lim_{d \to \infty} \frac{6|B_{2g}|[\NL_{g,d}^{\mathsf{prim}}]}{gd^{2g-1}\prod_{p \mid d}(1-p^{2-2g})} = \lambda_{g-1}\,.
    \end{equation}
    Denote the closures of $C$ and $C^{\mathsf{prim}}$ on $(\mathbf{NL}_g)_{\RR}$ by $\overline{C}$ and  $\overline{C^{\mathsf{prim}}}$. By \eqref{eq:limits}, we have $\lambda_{g-1} \in \overline{C}\subset \overline{C^{\mathsf{prim}}}$.
    
    Suppose that $\lambda_{g-1}$ lies in the boundary of $\overline{C}$. Since $\overline{C}$ is convex, there must be a nonzero linear functional $\ell : (\mathbf{NL}_g)_{\RR}\to \RR$ which vanishes on $\lambda_{g-1}$ and is nonnegative in $\overline{C}$, so $\ell (\Theta_{\NL}(q))$ is a cusp form all of whose Fourier coefficients are nonnegative, which is impossible unless $\ell ([\NL_{g,d}])=0$ for all $d$, a contradiction.
    
    Therefore, $\lambda_{g-1}$ lies in the interior of $\overline{C}$. Since $C$ is convex, it has the same interior as its closure, so $\lambda_{g-1}$ lies in the interior of $C$, and we can find some $N$ such that the cone generated by $[\NL_{g,1}], \ldots, [\NL_{g,N}]$ is full dimensional and $\lambda_{g-1}$ lies in its interior. Then, by \eqref{eq:limits}, there is some $N'>N$ such that $[\NL_{g,d}]$ is contained in the cone generated by $[\NL_{g,1}], \ldots , [\NL_{g,N}]$ for all $d > N'$. Therefore, $C$ agrees with the cone generated by $[\NL_{g,1}], \ldots , [\NL_{g,N'}]$, so it is a rational polyhedral cone. The same argument applies to $C^{\mathsf{prim}}$.
\end{proof}

\begin{proof}[Proof of Theorem \ref{thm:curves}]
    By \cite{GL24}, there is a Cartesian diagram
    $$
    \begin{tikzcd}
	{\mathcal H_{g \to 1, d}} & {\NL_{g,d}} \\
	{\mathcal M_g} & {\mathcal A_g}
	\arrow[from=1-1, to=1-2]
	\arrow["{\mathsf{ft}}"', from=1-1, to=2-1]
	\arrow["i", from=1-2, to=2-2]
	\arrow["\Tor"', from=2-1, to=2-2]
    \end{tikzcd}
    $$
    where $\Tor : \Mc_g \to \Ac_g$ is the Torelli map. Note that by Riemann--Hurwitz, $\mathcal H_{g \to 1, 1}$ is empty if $g \geq 2$, and that by \cite{Loo95}, the class $\lambda_{g-1}$ vanishes on $\mathcal M_g$. Therefore,
    $$
    \Tor^* (\Theta_{\NL}(q)) = \sum_{d \geq 2} [\mathcal H_{g \to 1, d}] q^d
    $$
    is a modular form by Theorem \ref{thm:main}, and it is clearly a cusp form.
\end{proof}
\begin{proof}[Proof of Corollary \ref{cor:curve-vanish}]
    Any modular form of weight $2g$ whose first two coefficients vanish has to be a multiple of the discriminant form $\Delta(q)^2$ by another modular form of weight $2g-24$. These do not exist if $2g \leq 22$ or $2g =26$. Therefore the power series from Theorem \ref{thm:curves} has to vanish.
\end{proof}

\begin{proof}[Proof of Corollary \ref{cor:Lian-smooth}]
    Any modular form of weight $24$ whose first two Fourier coefficients vanish is a scalar multiple of $\Delta(q)^2$. Therefore,
    $$
    \sum_{d \geq 2} [\mathcal H_{12\to 1,d}]q^d = \Delta(q)^2[\mathcal H_{12\to 1, 2}]
    $$
    and it was shown in \cite{ACCHLMT25} that $[\mathcal H_{12\to 1, 2}]$ (the bielliptic cycle) is not tautological.
\end{proof}

\bibliographystyle{alpha}
\bibliography{main_s0}

\end{document}